\documentclass[11pt]{article}
\usepackage{amsmath}
\usepackage{amssymb}
\usepackage{graphicx}
\usepackage{color}
\def\n{\noindent}
\def\L{{\bf L}}

\def\D{{\mathcal D}}

\def\dint{\int\!\!\int}
\def\ve{\varepsilon}
\def\n{\noindent}

\def\Tilde{\widetilde}

\def\E{{\mathcal E}}
\def\W{{\mathcal W}}
\def\C{{\mathcal C}}

\def\F{{\mathcal F}}

\def\implies{\Longrightarrow}
\def\ds{\displaystyle}
\def\sqr#1#2{\vbox{\hrule height .#2pt
\hbox{\vrule width .#2pt height #1pt \kern #1pt
\vrule width .#2pt}\hrule height .#2pt }}
\def\square{\sqr74}
\def\endproof{\hphantom{MM}\hfill\llap{$\square$}\goodbreak}
\def\nn{\nonumber}
\def\bega{\begin{array}}
\def\enda{\end{array}}
\def\begi{\begin{itemize}}
\def\endi{\end{itemize}}
\def\bes{\begin{eqnarray}}
\def\ees{\end{eqnarray}}
\def\O{{\cal O}}

\def\R{{\mathbb R}}

\def\ov{\overline}
\def\Tilde{\widetilde}
\def\forall{\hbox{for all }~}

\def\v{\vskip 1em}
\def\vs{\vskip 2em}
\def\vsk{\vskip 4em}

\def\T{{\mathcal T}}
\def\X{{\mathcal X}}

\def\be{\begin{equation}}
\def\beq{\begin{equation}}
\def\bel{\begin{equation}\label}
\def\eeq{\end{equation}}

\begin{document}
\title{\bf Lipschitz Metrics for a Class of Nonlinear Wave Equations}
\vs

\author{Alberto Bressan$^{(*)}$ and Geng Chen$^{(**)}$\\    \\
(*) Department of Mathematics, Penn State University,\\
University Park, Pa.~16802, U.S.A.\\ \, \\
(**) School of Mathematics
Georgia Institute of Technology,\\
Atlanta, Ga.~30332, U.S.A.
\\
e-mails:~ bressan@math.psu.edu,~gchen73@math.gatech.edu}

\maketitle
\begin{abstract} The nonlinear wave equation
$u_{tt}-c(u)(c(u)u_x)_x=0$ determines a flow of 
conservative solutions
taking values in the space $H^1(\R)$.    However, this flow is not continuous w.r.t.~the natural $H^1$ distance. Aim of this paper is to 
construct a new metric  which renders the 
flow uniformly Lipschitz continuous on bounded subsets of $H^1(\R)$.
For this purpose,  $H^1$ is given the structure of 
a Finsler manifold, where the norm of tangent vectors is defined
in terms of an optimal transportation problem.   For paths of 
piecewise smooth solutions,  one can carefully estimate how 
the weighted length grows in time.    By the generic regularity result
proved in \cite{BChen}, these piecewise regular 
paths are dense  and can be used to construct a
geodesic distance with the desired Lipschitz property. 
 \end{abstract}
\vsk

\section{Introduction} 
\setcounter{equation}{0}

Aim of this paper is to understand the continuous dependence
of solutions to the nonlinear wave equation
\bel{2} u_{tt}- c(u)\bigl(c(u) u_x\bigr)_x~=~0\,.\eeq
Roughly speaking, the analysis in \cite{BCZ, BZ, HR} shows that  conservative solutions are unique, globally defined, and
yield a flow on the space of couples
$(u,u_t)\in H^1(\R)\times \L^2(\R)$.  
For each conservative solution, the  total energy
\bel{en}
E(t)~\doteq~ \int\bigl[ u_t^2 + c^2(u)\, u_x^2\bigr]\, dx\eeq
remains constant in time. Precise results in this direction
will be recalled in Section~2.
On the other hand, these solutions
do not depend continuously on the initial data, w.r.t.~the
distance in the normed space $H^1\times \L^2$.


In the present paper we construct a new distance functional 
which renders  Lipschitz continuous the flow generated by (\ref{2}).
We recall that, for solutions of the Hunter-Saxton or
the Camassa-Holm equation, 
a similar task was achieved in \cite{BC1, BF, BHR, GHR1, GHR2}.

Developing ideas in \cite{BF}, our distance will be determined by the 
minimum cost to transport an energy measure from one solution to the other. While all previous papers dealt with first order equations,
to define a suitable transportation distance between
 two solutions $u, \tilde u$ of (\ref{2}) 
one now faces three main difficulties:
\begi
\item  At any given time $t$, each solution  determines two distinct
measures. These account for the energy $\mu^t_+$ 
of forward moving waves  and the 
energy $\mu_t^-$ of backward moving waves. 
The distance between $u(t)$ and $\tilde u(t)$ should be measured by 
the minimum  cost for  transporting $\mu^t_+$ to $\tilde \mu^t_+$
and $\mu^t_-$ to $\tilde \mu^t_-$.

\item The above double transportation problem is considerably complicated by
the fact that, while the total energy is conserved, some
energy can be transferred from forward to backward moving waves, or viceversa.
These source terms must be accounted for, when designing an ``optimal double
transportation plan".

\item As a wave front crosses waves of the opposite
family, its  speed can change.  As a consequence, the distance between two 
corresponding fronts in $u$ and $\tilde u$ may  quickly 
increase, making the optimal transportation plan more costly.  
To compensate for this effect, one needs to 
insert a weight function, accounting for the total energy of approaching waves.
\endi
In Section~3 we introduce
a Finsler norm on tangent vectors, related to an
energy transportation cost. Given a  smooth
path $\gamma:\theta\mapsto (u^\theta, u_t^\theta)$, 
one can then define
its weighted length $\|\gamma\|$ by integrating the norm of 
the tangent vector $d\gamma/d\theta$.   Proposition~1, stated in Section~3 and proved in Section~4,
contains the key 
estimate, describing how the norm of a tangent vector grows in time.
Assuming that, for $\theta\in [0,1]$ and 
$t\in [0,T]$, all solutions $u^\theta (t,\cdot)$
remain  sufficiently regular  so that the length 
of the path $\gamma^t:\theta \mapsto (u^\theta(t), u_t^\theta(t))$
 can still be computed, we obtain the bound
\bel{pl1}\|\gamma^t\|~\leq~C_T\,\|\gamma^0\|\,,\qquad\qquad\forall
t\in [0,T].
\eeq
Here the constant $C_T$ depends only on $T$ and on
a bound on the $H^1\times\L^2$ norm of the initial data.
At this stage, it is natural to define the geodesic distance
\bel{d1}
d^*\Big( (u, u_t)\,,~(\tilde u, \tilde u_t)\Big)~\doteq~
\inf~\Big\{ \|\gamma\|\,;~~~\gamma:[0,1]\mapsto H^1\times\L^2\,,
\quad \gamma(0)=(u, u_t),~~\gamma(1)=(\tilde u, \tilde u_t)\Big\}.\eeq
By (\ref{pl1}) we thus expect that, for any two solutions of 
(\ref{2}) and any $t\in [0,T]$, this distance  should satisfy
\bel{d2}d^*\Big( (u(t), u_t(t))\,,~(\tilde u(t), \tilde u_t(t))\Big)
~\leq~C_T\cdot 
d^*\Big( (u(0), u_t(0))\,,~(\tilde u(0), \tilde u_t(0))\Big).\eeq
This would imply that solutions depend Lipschitz continuously 
on the initial data, in the distance $d^*$.

\begin{figure}[htbp]
   \centering
\includegraphics[width=0.98\textwidth]{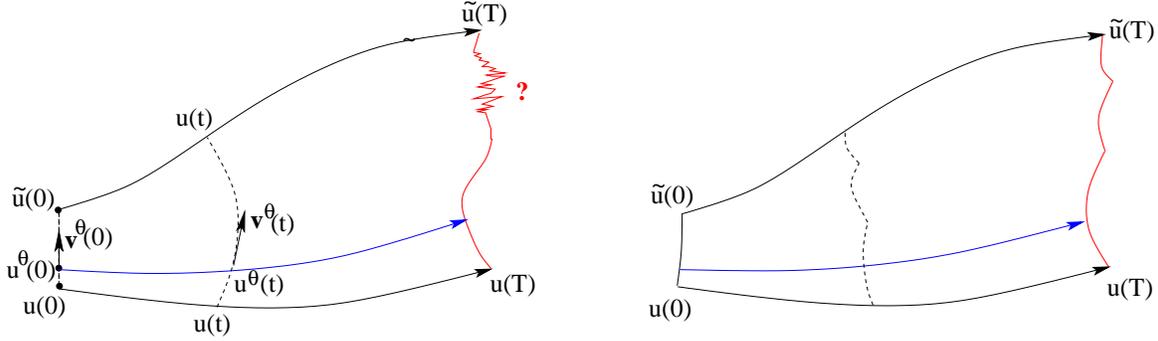}
   \caption{\small Left: due to singularity formation, a smooth path
  of initial data $\gamma^0:\theta\mapsto u^\theta(0)$ may lose regularity at a later time $T$.   In this case, the weighted length 
  $\|\gamma^T\|$ can no longer
  be computed by integrating the norm of a tangent vector.
  Right: by a small perturbation of the initial data, one obtains a path of solutions $\theta\mapsto u^\theta$ which remain piecewise smooth,
  for all except finitely many values of $\theta\in [0,1]$ .}
   \label{f:w26}
\end{figure}

To clinch this argument, one major difficulty 
must be overcome.  Indeed, smooth 
solutions may well develop singularities in finite time, \cite{GHZ}. Given a 
 path $\gamma^0$ 
 of smooth initial data, 
 there is no guarantee that at any time $t\in [0,T]$
 the path $\gamma^t$ will be regular enough so that the 
tangent vectors $d\gamma^t/d\theta$ are meaningfully defined (see Fig.~\ref{f:w26}).
We remark that a similar issue was encountered in the analysis
of hyperbolic conservation laws \cite{Bbook}.   For a path of piecewise
smooth solutions with finitely many shocks, a weighted norm 
on a suitable family of 
tangent vectors was introduced in \cite{BPisa}.
However, a lengthy effort was later required \cite{BCol, BCP}, 
in order to construct
paths of approximate solutions which retained enough regularity,
so that their length could still 
be estimated in terms of these tangent vectors.

In the present context, we can take advantage of the generic regularity results recently proved in \cite{BChen}. These can be summarized as follows. 
\begi
\item[(i)] For an open dense set of initial data 
\bel{idd}
(u(0,\cdot)\,, ~u_t(0,\cdot)) ~=~(u_0, u_1)~\in~
\Big(\C^3(\R)\cap H^1(\R)\Big) \times\Big(\C^2(\R)\cap\L^2(\R)\Big)
\eeq
the corresponding solution $u=u(t,x)$ of (\ref{2}) 
is piecewise smooth in the 
$t$-$x$ plane,  with 
singularities occurring along a finite set of smooth curves.
\item[(ii)] Every path of initial data $\theta\mapsto \gamma^0(\theta)
=(u_0^\theta, u_1^\theta)$ can be approximated by a second path
$\theta\mapsto \tilde \gamma^0(\theta)
=(\tilde u_0^\theta, \tilde u_1^\theta)$ such that,
for all but finitely many values of $\theta\in [0,1]$,
the corresponding solution $\tilde u^\theta$ remains piecewise smooth
on the domain $[0,T]\times \R$.
\endi
Using this dense set of piecewise regular paths, we can thus define
a geodesic distance on the space $H^1\times \L^2$, 
with the desired Lipschitz property.  Our main results are contained 
in 
\begi
\item Proposition~1, which establishes the basic estimate
(\ref{es0}) on the size of tangent vectors.
\item Theorem~5, providing the bound (\ref{lengg}) 
on how the length of a path of solutions can grow in time.
\item Theorem~7, showing that, by (\ref{LP7}),
the flow generated by the wave equation
(\ref{2}) is Lipschitz continuous w.r.t.~the
geodesic distance $d^*$. 
\endi
We remark that, for hyperbolic conservation laws, the 
distance constructed in \cite{BPisa, BCol, BCP} 
is equivalent to the $\L^1$ distance.  On the contrary,
our new metric is not equivalent to the norm 
distance on $H^1\times \L^2$. 
The completion of  $H^1\times \L^2$ 
w.r.t.~the geodesic distance includes a family of measures.  
This should not come as a surprise, since it was already 
observed in \cite{BZ, HR} that conservative solutions
can occasionally be measure-valued.

In Section~7 we compare
the geodesic distance (\ref{d1}) with more familiar  distances 
found in the literature.
In one direction, we show that
$$d^*\Big( (u_0, u_1)\,,~(\tilde u_0, \tilde u_1)\Big)
~\leq~C\cdot \Big( \|u_0-\tilde u_0\|_{H^1} + 
\|u_0-\tilde u_0\|_{W^{1.1}}+
\|u_1-\tilde u_1\|_{\L^2}+ \|u_1- \tilde u_1\|_{\L^1} \Big)\,,$$
for some constant $C$.
On the other hand, let $\mu$ and $\tilde\mu$ be the positive measures
having densities respectively
\bel{mmt}u_t^2 + c^2(u) u_x^2\qquad\hbox{and}\qquad 
\tilde u_t^2 + c^2(
\tilde u) \tilde u_x^2\eeq
w.r.t.~Lebesgue measure.
Then the geodesic distance $d^*$ dominates the Wasserstein
distance between the two measures.  Namely
\bel{wass}
\sup\,\left\{\bigg|\int f\,d\mu - \int f d\tilde\mu\bigg|\,;~~
\|f\|_{\C^1}\leq 1\right\}~\leq~
d^*\Big( (u, u_t)\,,~(\tilde u, \tilde u_t)\Big)\,.\eeq

All of the present analysis is concerned with conservative solutions
to (\ref{2}). For dissipative solutions, 
studied in \cite{BH, GHZ, ZZ1, ZZ2}, the continuous dependence
for general initial data
in $H^1\times \L^2$ remains an open question.  
For scalar conservation laws, 
an entirely different approach to  continuous dependence, 
relying on an $\L^2$ 
formulation, was
developed in \cite{BBL, B1, B2}.
\v
\section{Conservative solutions to the nonlinear 
wave equation} 
\setcounter{equation}{0}

In this section we review the main results in 
\cite{BChen, BCZ, BZ} on the Cauchy problem for the quasilinear
second order wave equation
\bel{0.1}
u_{tt} - c(u)\big(c(u) u_x\big)_x=0\,,
\eeq
with initial data
\bel{0.2}
u(0,x)=u_0(x)\,,\qquad
u_t(0,x)=u_1(x)\,.
\eeq
Here $c:\R\mapsto \R_+$ is
a smooth, uniformly positive function, such that
\bel{c0}
c(u)~\geq~c_0~>~0\,.\eeq
Consider the variables
\bel{0.3}
\left\{
\begin{array}{rcl}
R & \doteq  &u_t+c(u)u_x\,, \\
S & \doteq  &u_t-c(u)u_x\,,
\end{array} \right.
\eeq
so that
\bel{0.4}
u_t={R+S\over 2}\,,\qquad\qquad u_x={R-S\over 2c}\,.
\eeq
By (\ref{0.1}), the variables $R,S$ satisfy
\bel{0.5}
\left\{
\begin{array}{rcl}
R_t-cR_x &= & {c'\over 4c}(R^2-S^2), \\ [3mm]
S_t+cS_x &= & {c'\over 4c}(S^2-R^2).
\end{array} \right.
\eeq
Multiplying the first equation in (\ref{0.5}) by $R$ and the
second one by $S$, one obtains balance laws for $R^2$ and $S^2$, namely
\bel{0.6}\left\{
\begin{array}{rcl}
(R^2)_t - (cR^2)_x & = & {c'\over 2c}(R^2S - RS^2)\, , \\ [3mm]
(S^2)_t + (cS^2)_x & = & - {c'\over 2c}(R^2S -RS^2)\,.
\end{array}
\right.
\eeq
As a consequence, for smooth solutions the following quantity is
conserved:
\bel{0.7}
e~\doteq ~u_t^2+c^2u_x^2~=~{R^2+S^2\over 2}\,.
\eeq
We think of $R^2/2$ and $S^2/2$ as the energy of backward and 
forward moving waves, respectively.  These are not separately conserved. 
Indeed, by (\ref{0.6}) energy is transferred from forward to backward waves, and viceversa. The main results on the existence of solutions
to the Cauchy problem can be summarized as follows.
\v
\n{\bf Theorem 1.} {\it Let $c:\R\mapsto \R$ be
a smooth function satisfying (\ref{c0}).  Assume that the initial data
$u_0$ in (\ref{0.2}) is absolutely continuous, and that
$(u_0)_x\in\L^2\,$,   ~$u_1\in\L^2$.
Then the Cauchy problem (\ref{0.1})-(\ref{0.2}) admits a weak solution
$u=u(t,x)$, defined for all $(t,x)\in\R\times\R$.
In the $t$-$x$ plane, the function $u$ is locally H\"older continuous
with exponent $1/2$.   This solution $t\mapsto u(t,\cdot)$ is
continuously differentiable as a map with values in $\L^p_{\rm loc}$, for all
$1\leq p<2$. Moreover, it is Lipschitz continuous w.r.t.~the $\L^2$
distance, i.e.
\beq
\big\|u(t,\cdot)-u(s,\cdot)\big\|_{\L^2}~\leq~ L\,|t-s|\label{1.lip}
\eeq
for all $t,s\in\R$.
The equation (\ref{0.1}) is satisfied in distributional sense, i.e.
\beq
\dint \Big[\phi_t\, u_t - \big(c(u) \phi\big)_x c(u)\,u_x\Big]\,dxdt ~=~0
\label{consol}
\eeq
for all test functions $\phi\in\C^1_c$.
The maps $t\mapsto u_t(t, \cdot)$ and $t\mapsto u_x(t,\cdot)$ are continuous
with values in $\L^p_{loc}(\R)$, for every $p\in [1,2[\,$.}
\v
{\bf Theorem 2.} {\it In the same setting as Theorem~1,
 a unique  solution $u=u(t,x)$ exists
 which is {\em conservative} in the following sense.

 There exists two  families of positive Radon measures
 on the real line: $\{\mu_-^t\}$ and $\{\mu_+^t\}$, depending continuously
on $t$ in the weak topology of measures, with the following properties.
\begi
\item[(i)] At every time $t$ one has
\bel{E0}\mu_-^t(\R)+\mu_+^t(\R)~=~E_0~\doteq~
\int_{-\infty}^\infty \Big[u_1^2(x) + \bigl(c(u_0(x)) u_{0,x}(x)\bigr)^2\Big]\, dx \,.\eeq

\item[(ii)] For each $t$, the absolutely continuous parts of $\mu_-^t$ and
$\mu_+^t$ w.r.t.~the Lebesgue measure
have densities  respectively given  by
\bel{dmu}R^2 ~=~
\bigl(u_t + c(u) u_x\bigr)^2,\qquad\qquad
S^2~=~\bigl(u_t - c(u) u_x\bigr)^2.
\eeq
\item[(iii)]  For almost every $t\in\R$, the singular parts of $\mu^t_-$ and $\mu^t_+$
are concentrated on the set where $c'(u)=0$.
\item[(iv)] The measures $\mu_-^t$ and $\mu_+^t$ provide measure-valued solutions
respectively to the balance laws
\bel{mb4}\left\{
\begin{array}{rcl}
\xi_t - (c\xi)_x & = & {c'\over 2c}(R^2S - RS^2)\, , \\ [3mm]
\eta_t + (c\eta)_x & = & - {c'\over 2c}(R^2S -RS^2)\,.
\end{array}
\right.
\eeq
\endi
}
\v
The existence part of the above theorems was proved in \cite{BZ}.
The uniqueness of conservative solutions has been recently 
established in \cite{BCZ}.
\v
{\bf Remark 1.} By (\ref{mb4})
the total energy, represented by the positive measure 
$\mu^t=\mu^t_+ + \mu^t_-$, is
conserved in time.  Occasionally, some of this energy is concentrated
on a set of measure zero. At a time $\tau$ when this happens,
$\mu^\tau$ has a non-trivial singular part and
hence its absolutely continuous part satisfies
$$
\int \Big[u_t^2(\tau,x) +c^2\big(u(\tau,x)\big) \,u_x^2(\tau,x)\Big]
\,dx~<~E_0\,.$$
  The condition (iii) puts some restrictions
on the set of such times $\tau$. In particular, if $c'(u)\not= 0$
for all $u$, then this set has measure zero.
\v
{\bf Remark 2.} For any $t\geq 0$, the 
conservation of the total energy
implies
\bel{EE}\|u_t(t)\|^2_{\L^2}~
\leq~E_0~\doteq~\int (u_1^2 + c^2(u_0)u^2_{0,x})\, dx\,.\eeq
Hence  (\ref{1.lip}) holds with Lipschitz constant $L=\sqrt{E_0}$.
Moreover, one has the bounds
\bel{BB}
\|u(t,\cdot)\|_{\L^2}~\leq~\|u_0\|_{\L^2} + t\, \sqrt{E_0}\,,\qquad
\qquad \|u_x(t,\cdot)\|_{\L^2} ~\leq~{\sqrt{E_0}\over c_0}\,.\eeq
This yields an a priori bound on $\|u(t,\cdot)\|_{H^1}$, and hence on
$\|u(t,\cdot)\|_{L^\infty}$, depending only on time and on the total
energy $E_0$. In turn, since the wave speed $c(\cdot)$ is 
smooth, we obtain an a priori bound on $c(u)$ and $|c'(u)|$.
\v
\section{First order variations}
\setcounter{equation}{0}
For simplicity, in this section we consider solutions of (\ref{0.1}) with bounded support.   
More precisely, we shall assume that all our solutions satisfy
\bel{bsupp}u(t,x)~=~0 \qquad\qquad
\hbox{for}~ x\notin [0,L_0],\quad t\in [0,T].\eeq 
 Because of finite propagation
speed, this is hardly a restriction.

Let $(u, R,S)$ provide a smooth solution to (\ref{0.1}), (\ref{0.3}),  
and consider a family of perturbed
solutions of the form
\bel{perts}
u^\ve ~=~u+\ve v+o(\ve)\,,\qquad\qquad \left\{\bega{rl}R^\ve &=~R+\ve  r+o(\ve)\,,
\cr S^\ve &=~S+\ve   s + o(\ve)\,.\enda\right.\eeq
{}From (\ref{0.4}) it follows 
\bel{pvts}
u^\ve_t~=~{R^\ve+S^\ve\over 2}~=~{R+S\over 2} + \ve {  r+  s\over 2}+o(\ve)\,,
\eeq
\bel{prsx}
u^\ve_x~=~{R^\ve-S^\ve\over 2c(u^\ve)}~=~{R-S\over 2c(u)} + 
\ve {  r-  s\over 2c(u)}- \ve{R-S\over 2 c^2(u)} c'(u)\,v+o(\ve)\,.\eeq
Under the assumption (\ref{bsupp}), given $  r,   s$, the perturbation $v$ is uniquely
determined by 
\bel{ODEv}
v_x~=~-{(R-S)c'(u)\over 2c^2(u)}\,v +{  r-  s\over 2c(u)}\,,\qquad\qquad v(t,0)~=~0\,.\eeq
Furthermore, we have
\bel{ODEv2}
v_t~=~  {r+ s\over 2}\,.
\eeq
A direct computation shows that the first order perturbations $v,  s,  r$ satisfy the linear equations
\bel{le1}
v_{tt}- c^2 v_{xx} ~=~{2c c' u_x v_x} + \Big( (c')^2 u_x^2 + c c''u_x^2 + 2 c c' u_{xx}
\Big) v\,.\eeq
\bel{le2}\left\{\bega{rl}   r_t - c(u)   r_x&\ds =~c' R_x v + \left( 
{c''\over 4c}- {(c')^2 \over 4c^2}\right) (R^2-S^2) v + {c'\over 2c} 
(R  r-S  s)\,,\cr
&\cr
   s_t + c(u) s_x&\ds=~-c' S_x v + \left( 
{c''\over 4c}- {(c')^2 \over 4c^2}\right) (S^2-R^2) v + {c'\over 2c}
 (S  s-R  r)\,.\enda\right.
\eeq
By the assumptions (\ref{c0}) on 
 the wave speed  $c(u)$,  all functions
$c'/4c$, $c''/ 4c$, $(c')^2/ 4c^2$, are smooth functions of $u$. 
\v
We shall introduce a weighted norm on tangent vectors $r,s$,  which takes into account 
the total  energy of waves which are approaching
a given wave located at $x$.  This is described by the weights
\bel{we} \W^-(x) ~\doteq~\ds 1+\int^x_{-\infty} S^2(y)\, dy\,,\qquad\qquad
\W^+(x) ~\doteq~\ds 1+\int_x^{+\infty} R^2(y)\, dy
\,.\eeq
In addition, consider the function
\bel{adef}a(t)~\doteq~\int_{-\infty}^\infty\frac{|c'|\,\bigl|R^2S-S^2R
\bigr|}{2c}(t,x)~dx.\eeq
As proved in \cite{BCZ},
the function 
$$\tau~\mapsto~\int_0^\tau \int_{-\infty}^{+\infty}
 \left|\frac{c'}{2c}(R^2S-RS^2)\right|(t,x)
\,dx\, dt$$
is  H\"older continuous and absolutely continuous on bounded 
time intervals, and has sub-linear growth.
In particular (see (3.11)-(3.12) in the proof of Lemma~1
in \cite{BCZ}), one has 
\bel{iabd}\int_0^T a(t)\, dt ~\leq~C_T\,,\eeq
for some constant $C_T$ depending only on $T$  and on the
total energy $E_0$.
By (\ref{0.6}) it follows
\bel{Wd}\left\{\bega{rl}
\W^-_t-c\W^-_x&\ds=~ -2cS^2+\int^x_{-\infty}
\frac{c'}{2c}(S^2R-R^2S)~dy~\leq~-2c_0S^2+ a(t),\cr\cr
\W^+_t+c\W^+_x&\ds=~-2cR^2
+\int^{+\infty}_{x}\frac{c'}{2c}(R^2S-S^2R)~dy~\leq~-2c_0R^2+a(t)\,.
\enda\right.\eeq

On the space of tangent vectors $(v,r,s)$ we introduce a 
Finsler  norm, having the form 
\bel{fn1}
\Big\|(v,\,r,\,s)\Big\|_{(u,R,S)}
~\doteq~\inf_{\tilde r,\tilde s, w,z}
\Big\|(\tilde r,w,\tilde s,z)\Big\|_{(u,R,S)}\,,
\eeq
where the infimum is taken over the set of 
vertical displacements $\tilde r,\tilde s$ and shifts $w,z$
which satisfy
\bel{rstt}\left\{ 
\bega{rl} r&=~\tilde r - w R_x + {c'\over 8c^2}(w-z)S^2\,,\cr\cr
s&=~\tilde s - z S_x + {c'\over 8c^2} (w-z)R^2\,.\enda
\right.\eeq
This norm is defined as
\bel{main}\bega{l}
\Big\|(\tilde{r},w,\tilde{s},z) \Big\|_{(u,R,S)}\cr\cr
\doteq
\ds~ \kappa_1 \int \Big\{   |w|\bigl(1+R^2\bigr)\W^-  +  
|z|\bigl(1+S^2\bigr)\W^+    \Big\}dx\cr\cr
\quad\ds+\kappa_2\int   \Big\{   |\tilde{r}|\W^-  + 
 |\tilde{s}|\W^+    \Big\}dx\cr\cr
\quad \ds
+\kappa_3\int  \Big|v+\frac{Rw}{2c}-\frac{Sz}{2c}\Big|\,\bigg\{(1+R^2)\,\W^- 
+
(1+S^2)\,\W^+ \bigg\}dx 
\cr\cr\quad\ds+\kappa_4 \int   \bigg\{   \Big|w_x+\frac{c'}{4c^2}(w-z)S
\Big|\W^-  +  \Big|z_x+\frac{c'}{4c^2}(w-z)R\Big|
\W^+    \bigg\}dx
\cr\cr
\quad\ds
+\kappa_5\int   \bigg\{   \Big|Rw_x+\frac{c'}{4c^2}(w-z)SR
\Big|\W^- 
 +  
\Big|Sz_x+\frac{c'}{4c^2}(w-z)RS\Big|\W^+    \bigg\}dx
\cr\cr
\ds\quad+\kappa_6\int   \bigg\{ \Big|2R\tilde{r}+R^2w_x+\frac{c'}{4c^2}R^2S(w-z)
\Big|\W^- 
+   \Big|2S\tilde{s}+S^2z_x+\frac{c'}{4c^2}S^2R(w-z)\Big|\W^+    
\bigg\}dx
\cr\cr
\doteq~\kappa_1 I_1+\kappa_2I_2+\kappa_3I_3+\kappa_4I_4+\kappa_5I_5+ \kappa_6I_6\,,\enda\eeq
for suitable constants $\kappa_1,\ldots,\kappa_6$ to be determined later.

\begin{figure}[htbp]
   \centering
\includegraphics[width=0.5\textwidth]{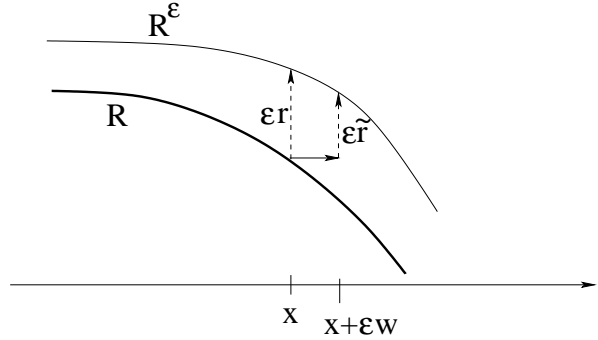}
   \caption{\small A perturbation of the $R$-component of the solution to the
   variational wave equation.}
   \label{f:hyp37}
\end{figure}

 To motivate (\ref{fn1}), consider
a profile $R$ and a perturbation $R^\ve$, as shown in figure~\ref{f:hyp37}.
In first approximation, $R^\ve\approx R+\ve r$.    Notice that we could also
obtain the profile $R^\ve$  starting from 
the graph of $R$, performing a horizontal shift in the amount $\ve w$ and  then 
a vertical shift in the amount $\ve \tilde r$, provided that
\bel{tvr}
r ~= ~\tilde r - w R_x\,.\eeq
As a first guess, one could thus 
define a norm $\|r\|_\dagger $ by optimizing the choice
of $\tilde r, w$, subject to (\ref{tvr}).   However, a detailed analysis has shown
that this approach does not  work. Indeed, it does not take into account
the fact that, when  backward and  forward moving waves 
cross each other,  by (\ref{0.5})
their sizes $R,S$ are modified.
Compared with (\ref{tvr}), the additional term 
in the first equation of (\ref{rstt})
accounts for this interaction.    Notice that 
$w-z$ is the relative shift of backward w.r.t.~forward waves.

We now explain the meaning of each integral on the right hand side of
(\ref{main}).   
\begi
\item  The integral of  $|w|(1+R^2)$ can be interpreted as the  
cost for  transporting the base measure
with density  $1+R^2$  
 from the point $x$ to the point $x+\ve w(x)$. 

Similarly, the integral of $|z|(1+S^2)$ accounts for the cost of transporting
the measure with density $1+S^2$ from $x$ to $x+\ve z(x)$. 

Here, as in all other terms, we  insert the weights $\W^\pm$ 
coming from the interaction potential.
\v

\item $I_2$ accounts for the vertical shifts in the graphs of $R,S$.
We interpret the integrand as the change in $\arctan R$ times the 
density $(1+R^2)$ of the base measure.     Notice that here the factor $(1+R^2)$ 
cancels out with the derivative of the arctangent.

\item $I_3$ accounts for the changes in $u$. Observe that
$$\ve^{-1}[u^\ve(x+\ve w(x))- u(x)]~\approx~ v(x) +  u_x(x) w(x)~=~
v(x) +  {R(x)-S(x)\over 2c(u(x))} w(x). $$
This can be written in the form
\bel{q1}v +  {R-S\over 2c} w~=~\left(v +  {R\over 2c} w - {S\over 2c} z
\right) +{S (z -w)
\over 2c}\,.\eeq
Notice that the last term on the right hand side of (\ref{q1}) does not appear  
in $I_3$. In fact, the last term  ${S(z-w)\over 2c}$ is the relative shift term coming from the equation \eqref{0.3}. Subsequent computations will show that 
this term is inessential,
because its contribution
can be bounded by the decrease in the interaction potential.  In an entirely similar way
we obtain
$$\ve^{-1}[u^\ve(x+\ve z(x))- u(x)]~\approx~ v(x) +  u_x(x) z(x)~=~
v(x) +  {R(x)-S(x)\over 2c(u(x))} z(x), $$
$$v +  {R-S\over 2c} z~=~\left(v +  {R\over 2c} w - {S\over 2c} z
\right) +{R (z -w)
\over 2c}\,.$$

\item $I_6$ accounts for the change in base measure with densities
$R^2$ and $S^2$,
produced by the shifts $w,z$. To see this, assume that the mass with density 
$R^2$ is transported from $x$ to $x+\ve w(x)$.   If the mass were conserved, 
the new density should be
\bel{t1}(R^\ve)^2(x+\ve w(x))~=~R^2(x) -\ve w(x) 
(R^2)_x(x) - \ve w_x(x) R^2(x)+o(\ve).
\eeq
In addition, if the mass with density $S^2$ is transported
from $x$ to $x+\ve z(x)$, by (\ref{0.6})
the crossing between forward and backward waves yields
the source term
\bel{t2}{c'\over 2c}(R^2S-RS^2)\cdot {z-w\over 2c}\,.\eeq
On the other hand, if we shift the graph of $R$ horizontally by $\ve w$ 
and then vertically by $\ve \tilde r$, the 
new density will be
\bel{t3}(R^\ve)^2(x+\ve w(x))~=~R^2(x) -\ve w(x) (R^2)_x(x) + 2\ve R(x)\tilde r(x)+o(\ve).\eeq
Subtracting (\ref{t1})-(\ref{t2}) from  (\ref{t3}) we obtain the expression
\bel{t4}2R( r+w R_x ) + R^2 w_x +{c'\over 4c^2}(R^2S-RS^2)\cdot (w-z)\,.\eeq
\item The integrals $I_4$ and $I_5$
does not seem to have a clear geometric 
interpretation. $I_4$ is somewhat related to  
the change in Lebesgue measure
produced by the shifts $w,z$, while
$I_5$ is related to the change in base measure with densities
$R$ and $S$,
produced by the shifts $w,z$.   As shown by our
subsequent computations, these two additional terms 
must be included in the 
definition (\ref{main}), in order to estimate the 
 time derivatives of  $I_3$ and $I_6$.
\endi
\v
Our goal is to prove

{\bf Proposition 1.}  {\it Let $(u,R,S)$ be a smooth solution to 
(\ref{0.1}) and (\ref{0.5}), and assume that the first order perturbations
$(v,r,s)$ satisfy the corresponding 
linear equations (\ref{le1})-(\ref{le2}).
Then for any $\tau\geq 0$ one has
\bel{es0}
 \Big\|(v(\tau), r(\tau), s(\tau))\Big\|_{(u(\tau),R(\tau),S(\tau))}~\leq~
\exp\left\{C\tau+\int_0^\tau a(s)ds\right\}
\cdot\Big\|(v(0), r(0), s(0))\Big\|_{(u(0),R(0),S(0))}\,,\eeq
 with
a constant $C$ depending only on the total energy.}
\v
Toward the proof, 
the main argument goes as follows.
At time $t=0$ let a tangent vector 
$
(v(0),r(0), s(0))$
be given.  By the definition (\ref{fn1}), for any $\epsilon>0$
we can find shifts  $w_0,z_0$ and perturbations $\tilde r_0,\tilde s_0$
satisfying
\bel{fn6}\Big\|(\tilde r_0,w_0,\tilde s_0,z_0)\Big\|_{(u(0),R(0),S(0))}~\leq~
\epsilon+
\Big\|(v(0),\,r(0),\,s(0))\Big\|_{(u(0),R(0),S(0))}
\eeq
together with the constraints
\bel{rstt3}\left\{ 
\bega{rl} r(0)&=~\tilde r_0 - w_0 R_x(0) + {c'\over 8c^2}(w_0-z_0)S^2(0)\,,\cr\cr
s(0)&=~\tilde s_0 - z_0 S_x(0) + {c'\over 8c^2} (w_0-z_0)R^2(0)\,.\enda
\right.\eeq
In order to prove (\ref{es0}), for any $t\in [0,\tau]$ it suffices to find
shifts $w(t), z(t)$, together with $\tilde r(t)$, $\tilde s(t)$ satisfying
(\ref{rstt}) and the initial condition (\ref{rstt3}),  so that
\bel{wn1}
{d\over dt} \Big\|(\tilde r(t), w(t),\tilde s(t), z(t))\Big\|_{(u(t),R(t),S(t))}~\leq~
\bigl(C+ a(t)\bigr)
\cdot\Big\|(\tilde r(t), w(t),\tilde s(t), z(t))\Big\|_{(u(t),R(t),S(t))} \,.\eeq
These shifts $w(t), z(t)$ will be obtained by 
propagating  along characteristics 
the shifts $w_0,z_0$ in the initial data. More precisely, we choose
$w,z$ to be the solutions of the linearized system
\bel{wzt}\left\{
\bega{rl} w_t - c(u) w_x &=~-c'(u)\,(v+u_xw)\,,\cr\cr
z_t + c(u) z_x&=~c'(u) (v+u_x z)\,,\enda\right.\eeq
with initial data
\bel{wzid}\left\{
\bega{rl} w(0,x)&=~ w_0(x)\,,\cr\cr
z(0,x)&=~ z_0(x)\,.\enda\right.\eeq
By (\ref{le2}) and the identities (\ref{rstt}), this determines
the evolution equation for $\tilde r, \tilde s$.

In the next section, by carefully estimating the time derivatives 
of all terms in
(\ref{main}), we shall prove that (\ref{wn1})  holds.
In turn, this will yield (\ref{es0}).
\v
\section{Estimates on the norm of tangent vectors\label{section_3}}
\setcounter{equation}{0}

The first part of the proof of (\ref{wn1}) is largely computational.
Using the evolution equations (\ref{0.1}), (\ref{0.3}), (\ref{0.5})
for $u,R,S$, and (\ref{le2}), (\ref{wzt}) for $r,s,w,z$, together with the 
identities (\ref{rstt}), we estimate the time derivative of each
integral in (\ref{main}).
\v
{\bf 1.}
To estimate the time derivative of $I_1$ (shift in the base measure), 
using (\ref{wzt}) 
we first compute
$$\bega{l}
\bigl( w(1+R^2) \bigr)_t -\bigl( c w(1+R^2) \bigr)_x\cr\cr
\qquad =~(w_t-c w_x)(1+R^2) + w\bigl[  (R^2)_t -(cR^2)_x  \bigr] - w c_x\cr\cr
\qquad \ds=~-c'\Big(v+\frac{R-S}{2c} w\Big)(1+R^2) + \frac{c'}{2c}w(R^2 S- R S^2- R+ 
S)\cr\cr
\qquad \ds=~-c'
\Big(v+\frac{R}{2c} w-\frac{S}{2c}z\Big)(1+ R^2) + \frac{c'}{2c}w(2R^2 S- R S^2- R+2 
S)
-\frac{c'}{2c}zS(1+R^2).
\enda
$$
Thanks to (\ref{Wd}) we obtain 
\bel{DI1}
\bega{l}\ds\frac{d}{dt}\int    |w|\bigl(1+R^2\bigr)\W^- \, dx
~\leq~ O(1)\cdot\int |w|\,\bigl(1+|R^2 S| + |RS^2| +|R|+|S|\bigr)\W^-\, dx \cr\cr
\qquad \ds+
O(1)\cdot\int |z|\,\bigl(|S| +  |R^2S|\bigr)\W^+\, dx+ \O(1)\cdot \int 
\Big|v+\frac{Rw}{2c} -{Sz\over 2c}\Big|(1+R^2)\W^-\, dx
\cr\cr
\ds \qquad  +a(t)\int    |w|\bigl(1+R^2\bigr)\,\W^- \, dx-
 2c_0\int   |w|\bigl(1+R^2\bigr) S^2\W^- \, dx\,.\enda
\eeq
\v
{\bf 2.} To estimate the time derivative of $I_2$ (change in arctan), using 
(\ref{le2}) we first compute
\bel{3.1}\bega{l}
\bigl( r+w R_x \bigr)_t -\bigl( c(r+w R_x) \bigr)_x\cr\cr
\qquad \ds=~\bigl[r_t-(c r)_x\bigr]
 + (w_t -c w_x) R_x+ w\bigl[  (R_x)_t -(cR_x)_x  \bigr]\cr\cr
\qquad \ds = ~- c'\,{R-S\over 2c}\,r+c' R_x v +\frac{c''c-(c')^2}{4c^2}(R^2-S^2)v 
+{c'\over 2c} (Rr-Ss)\cr\cr
\qquad\qquad \ds -c'\,\left( v+{R-S\over 2c}w\right)R_x 
+w\left[ \frac{c''c-(c')^2}{4c^2}\,\frac{R-S}{2c}\,(R^2-S^2)  +\frac{c'}{4c}(2RR_x-2SS_x) 
\right]
\cr\cr
\qquad =~\ds \frac{c'}{2c}Sw(R_x- S_x) + \frac{c'}{2c}S (r-s)
+\frac{c''c-(c')^2}{4c^2}(R^2-S^2)\left( v+ {R-S\over 2c}w\right)\,.\enda
\eeq
Next,
\bel{3.2}\bega{l}
\ds \left(\frac{c'}{8c^2}(w-z)S^2\right)_t-\left (c\frac{c'}{8c^2}(w-z)S^2\right)_x\cr\cr
\qquad \ds = ~\frac{c''c-2(c')^2}{8c^3}(u_t - c u_x) (w-z)S^2+\frac{c'}{8c^2}(w_t-cw_x)S^2 
-\frac{c'}{8c^2}(z_t+cz_x)S^2\cr\cr
\qquad\qquad \ds+\frac{c'}{8c^2}2cz_xS^2+
	\frac{c'}{8c^2}(w-z)\left[ (S^2)_t + (cS^2)_x \right]-\frac{c'}{8c^2}2(w-z)(cS^2)_x
	\cr\cr
\qquad \ds =~\frac{c''c-2(c')^2}{8c^3}(w-z)S^3-\frac{(c')^2}{8c^2}\left(v+\frac{R-S}{2c}w
\right)S^2-\frac{(c')^2}{8c^2}\left(v+\frac{R-S}{2c}z
\right)S^2
\cr\cr
\qquad\quad\  \ds
	+\frac{c'}{4c}z_xS^2-\frac{(c')^2}{16c^3}(w-z)(R^2 S - R S^2) 
	-\frac{(c')^2}{8c^3}(w-z)(RS^2-S^3)-\frac{c'}{2c}(w-z)SS_x
\qquad \ds 
	\,.\enda
\eeq
By (\ref{rstt}), combining (\ref{3.1}) with (\ref{3.2}) we obtain
\bel{3.3}\bega{l}\ds
\qquad\tilde r_t - (c \tilde r)_x
\cr
\cr
\ds=~\left[\bigl( r+w R_x \bigr)_t -\bigl( c(r+w R_x) \bigr)_x\,
\right]
-\left[\left(\frac{c'}{8c^2}(w-z)S^2\right)_t-\left (c\frac{c'}{8c^2}(w-z)S^2\right)_x\,
\right]
\cr
\cr
\ds 
=~\frac{c'}{2c}Sw(R_x- S_x) + \frac{c'}{2c}S (r-s)
+\frac{c''c-(c')^2}{4c^2}(R^2-S^2)\left( v+ {R-S\over 2c}w\right)\cr\cr
\qquad
\qquad\ds 
-\frac{c''c-2(c')^2}{8c^3}(w-z)S^3+\frac{(c')^2}{8c^2}\left(2v+\frac{R-S}{2c}(w+z)
\right)S^2-\frac{c'}{4c}z_xS^2
\cr
\cr
\qquad\qquad \ds
	+\frac{(c')^2}{16c^3}(w-z)(R^2 S - R S^2) 
	+\frac{(c')^2}{8c^3}(w-z)(RS^2-S^3)+\frac{c'}{2c}(w-z)SS_x
\cr\cr
 \ds = ~{c'\over 2c} (wSR_x-zSS_x)-\frac{c'}{4c}z_xS^2+ \frac{c'}{2c}S 
\left( (\tilde r-\tilde s) - (wR_x - z S_x) +{c'\over 8c} (w-z)(S^2-R^2)\right)\cr\cr
\qquad\qquad \ds+\frac{c''c-(c')^2}{4c^2}(R^2-S^2)\left( v+ {R-S\over 2c}w\right)
-\frac{c''c-2(c')^2}{8c^3}(w-z)S^3\cr\cr
\qquad\qquad
\ds 
+\frac{(c')^2}{8c^2}\left(2v+\frac{R-S}{2c}(w+z)
\right)S^2
\cr\cr
\ds=~{c'\over 2c} S \tilde r - {c'\over 4c}\left( 2S\tilde s + S^2 z_x + {c'\over 4c} 
S^2R(w-z)\right)\cr\cr
\ds\qquad\qquad +\frac{c''c-(c')^2}{4c^2} \, R^2\left( v+{Rw\over 2c}  - {Sz\over 2c}
\right) - 
\frac{c''c-2(c')^2}{4c^2} \, S^2\left( v+{Rw\over 2c}  - {Sz\over 2c}\right)\cr\cr
\qquad\qquad + \O(1)\cdot \bigl( |w| + |z|\bigr) \bigl(1+ |R^2 S|+ |R S^2|\bigr)\,.
\enda
\eeq
We thus conclude
\bel{DI2}\bega{l}
\ds\frac{d}{dt}\int |\tilde{r}|\W^- dx~=~
\O(1)\cdot \int |S\tilde r|\,\W^- + \O(1)\cdot \int \,\left| 2S\tilde s + S^2 z_x + {c'\over 4c} 
S^2R(w-z)\right|\, \W^+\, dx \cr\cr
\quad \ds +
\O(1)\cdot \int \,S^2\left| v+{Rw\over 2c}  - {Sz\over 2c}\right| \W^+\, dx+ 
\O(1)\cdot \int \,R^2\left| v+{Rw\over 2c}  - {Sz\over 2c}\right| \W^-\, dx
\cr\cr
\ds\quad
\ds+
\O(1)\cdot \int |w|  \bigl( 1+|R^2 S|+ |R S^2|\bigr)\,\W^-\, dx
  + 
\O(1)\cdot \int  |z| \bigl(1+ |R^2 S|+ |R S^2|\bigr)\,\W^+\, dx\cr\cr
\quad\ds +a(t)\int |\tilde{r}|\W^- dx- 2 c_0 \int  |\tilde r| \,S^2\,\W^-dx\,.
\enda
\eeq

{\bf 3.}
To estimate the time derivative of $I_3$ (change in $u$), using the identities 
in (\ref{ODEv})-(\ref{ODEv2}) for $v_t$ and $v_x$,
we first compute
\bel{3.4}
v_t-cv_x~= ~s +\frac{c'}{2c}(R-S) v.
\eeq
Next, by (\ref{0.3}) and (\ref{wzt})  we obtain
\bel{3.5}
\bega{l}\ds
\left(\frac{Rw}{2c}-\frac{Sz}{2c}\right)_t-c\left(\frac{Rw}{2c}-\frac{Sz}{2c}\right)_x
\cr\cr
\ds\qquad =~\frac{1}{2c}w(R_t-cR_x)-\frac{1}{2c}z(S_t+cS_x)+zS_x
\cr\cr
\ds\qquad\qquad 
+\frac{R}{2c}(w_t-cw_x)-\frac{S}{2c}(z_t+cz_x)+Sz_x-\frac{c'}{2c^2}(RSw-S^2z)
\cr\cr
\ds \qquad =~\frac{c'}{8c^2}w(R^2-S^2)-\frac{c'}{8c^2}z(S^2-R^2)+zS_x
\cr\cr
\ds
\qquad \qquad
-\frac{c'}{2c}R\Big(v+\frac{R-S}{2c}w\Big)
-\frac{c'}{2c}S\Big(v+\frac{R-S}{2c}z\Big)+Sz_x  -\frac{c'}{2c^2}(RSw-S^2z),\enda
\eeq
Finally, by (\ref{0.5}) it follows
\bel{3.6}
(1+R^2)_t-(c(1+R^2))_x~=~\frac{c'}{2c}(R^2S-RS^2)-\frac{c'}{2c}(R-S).\nn
\eeq
Putting together (\ref{3.4})--(\ref{3.6}) and using (\ref{rstt}) one obtains
\bel{3.7}\bega{l}
\ds\left[ \Big(v+\frac{Rw}{2c}-\frac{Sz}{2c}\Big)(1+R^2)\right]_t -\left[ c\Big(v+\frac{Rw}{2c}-\frac{Sz}{2c}\Big)(1+R^2)\right]_x
\cr\cr
\ds\qquad =~ \left[v_t-cv_x+\Big(\frac{Rw}{2c}-\frac{Sz}{2c}\Big)_t-
c\Big(\frac{Rw}{2c}-\frac{Sz}{2c}\Big)_x\right](1+R^2) 
\cr\cr
\qquad \qquad \ds+\Big(v+\frac{Rw}{2c}-\frac{Sz}{2c}\Big)
\Big[(1+R^2)_t-\bigl(c(1+R^2)\bigr)_x\Big]
\cr\cr
\qquad =\ds~\bigg[s +\frac{c'}{2c}v(R-S)+\frac{c'}{8c^2}w(R^2-S^2)
-\frac{c'}{8c^2}z(S^2-R^2)+zS_x\cr\cr
\qquad\qquad\ds
-\frac{c'}{2c}R\Big(v+\frac{R-S}{2c}w\Big)-\frac{c'}{2c}S\Big(v+\frac{R-S}{2c}z\Big)+Sz_x  -\frac{c'}{2c^2}(SRw-S^2z)\bigg](1+R^2)\cr\cr
\qquad\qquad\ds+\Big(v+\frac{Rw}{2c}-\frac{Sz}{2c}\Big)
\left[\frac{c'}{2c}(R^2S-RS^2)-\frac{c'}{2c}(R-S)\right]\cr\cr
\qquad\ds=~
\left[\tilde{s}+\frac{c'}{8c^2}(z-w)S^2-\frac{c'}{c}S\Big(v+\frac{Rw}{2c}-\frac{Sz}{2c}\Big)
+\Big(Sz_x+\frac{c'}{4c^2}(w-z)RS\Big)\right](1+R^2)
\cr\cr
\qquad\qquad\ds
+\Big(v+\frac{Rw}{2c}-\frac{Sz}{2c}\Big)
\left[\frac{c'}{2c}(R^2S-RS^2)-\frac{c'}{2c}(R-S)\right].
\enda
\eeq
We thus conclude
\bel{DI3}\bega{l}\ds
\frac{d}{dt}\int  \left|v+\frac{Rw}{2c}-\frac{Sz}{2c}\right|\,(1+R^2)\,\W^- \,
dx\cr\cr
\ds
\quad
\leq~ \int |\tilde s|(1+R^2)\W^+\, dx + 
\ds \O(1)\cdot\int  \left|v+\frac{Rw}{2c}-\frac{Sz}{2c}\right|\,
\bigl(1+|R|+|S|+|R^2S| + |RS^2|\bigr)\,\W^- 
\,dx 
\cr\cr
\ds\quad\qquad 
+\O(1)\cdot\int |w| S^2(1+R^2)\W^-\, dx
+\O(1)\cdot\int |z| S^2(1+R^2)\W^+\, dx\cr\cr
\qquad\qquad \ds 
+\O(1)\cdot\int \left|Sz_x+\frac{c'}{4c^2}(w-z)RS\right|\,(1+R^2)\W^-\, dx
\cr\cr
\quad\qquad \ds+a(t)\int  \left|v+\frac{Rw}{2c}-\frac{Sz}{2c}\right|\,(1+R^2)\,\W^- 
-2 c_0\int \left|v+\frac{Rw}{2c}-\frac{Sz}{2c}\right|\,(1+R^2)S^2\,\W^-dx\,.  \enda
\eeq

\vs

{\bf 4.}
To estimate the time derivative of $I_4$, recalling (\ref{wzt})
we first compute
\bel{4.1}\bega{l}
(w_x)_t-(cw_x)_x\cr\cr
\qquad \ds=~-\frac{c''}{2c}(R-S)\Big(v+\frac{R-S}{2c}w\Big)-c'\left[
-\frac{(R-S)c'}{2c^2}v+{ \frac{r-s}{2c}}
-\frac{c'}{4c^3}(R-S)^2 w \right]\cr\cr
\qquad \ds\qquad -\frac{c'}{2c}(R_x w-S_x w)-\frac{c'}{2c}(R-S)w_x\,.\enda
\eeq
Moreover, by (\ref{0.3}) and (\ref{wzt}), one has
\bel{4.2}\bega{l}\ds
\left(\frac{c'}{4c^2}wS\right)_t-\left(c\frac{c'}{4c^2}wS\right)_x\cr\cr
\quad \ds=~\left(\frac{c'}{4c^2}\right)'wS^2
-\frac{(c')^2}{4c^2}\left(v+\frac{R-S}{2c}w\right)S -\frac{(c')^2}{16c^3}w(R^2-S^2)-\frac{(c')^2}{8c^3}(R-S)wS-\frac{c'}{2c}wS_x\,,\enda
\eeq
\bel{4.3}
\bega{rl}\ds\left(\frac{c'}{4c^2}zS\right)_t-\left(c\frac{c'}{4c^2}zS\right)_x
&\ds=~\left(\frac{c'}{4c^2}\right)'zS^2+\frac{(c')^2}{4c^2}\left(v+\frac{R-S}{2c}z\right)
S +\frac{(c')^2}{16c^3}z(S^2-R^2)\cr\cr
& \qquad \ds-\frac{(c')^2}{8c^3}(R-S)zS-\frac{c'}{2c}zS_x-\frac{c'}{2c}z_x S\,.
\enda 
\eeq
Combining the identities (\ref{4.1})--(\ref{4.3}) and recalling (\ref{rstt}), we obtain
\bel{4.4}\bega{l}
\ds
\Big(w_x + {c'\over 4c^2}(w-z)S\Big)_t -\left[c\Big(w_x + {c'\over 4c^2}(w-z)S\Big)
\right]_t\cr\cr
\qquad
=~
\ds\frac{c'}{2c}\tilde{s}-\frac{c'}{2c}\tilde{r}+\frac{c'}{2c}S\Big(w_x+\frac{c'}{4c^2}(w-z)S\Big)\cr\cr
\qquad\qquad \ds
+\frac{c'}{2c}\Big(Sz_x+\frac{c'}{4c^2}(w-z)SR\Big)
-\frac{c'}{2c}\Big(Rw_x+\frac{c'}{4c^2}(w-z)SR\Big)
\cr\cr
\qquad\qquad \ds
-\frac{c''c-(c')^2}{2c^2}R\Big(v+\frac{Rw}{2c}-\frac{Sz}{2c}\Big)+
\frac{c''c-2(c')^2}{2c^2}S\Big(v+\frac{Rw}{2c}-\frac{Sz}{2c}\Big)
\cr\cr
\qquad\qquad \ds +\frac{c''c-(c')^2}{4c^3}RS(w-z)-\frac{(c')^2}{8c^3}S^2(w-z).\enda
\eeq
By the previous analysis, thanks to the uniform bounds (\ref{Wd})
on the weights, we conclude
\bel{DI4}\bega{l}\ds
\frac{d}{dt}\int \left|w_x+\frac{c'}{4c^2}(w-z)S\right|\W^-\, dx \cr\cr
\qquad\ds\leq~ O(1)\cdot \int |\tilde r|\,\W^-dx+O(1)\cdot \int |\tilde s|\,\W^+dx
+\O(1)\cdot\int \Big |S \Big | \Big |w_x+\frac{c'}{4c^2}(w-z)S\Big|\, \W^-dx
\cr\cr
\qquad\qquad \ds +\O(1)\cdot\int \Big|Sz_x+\frac{c'}{4c^2}(w-z)RS\Big|\,\W^+dx 
+\O(1)\cdot\int \Big|Rw_x+\frac{c'}{4c^2}(w-z)RS\Big|\, \W^-dx
\cr\cr
\qquad\qquad \ds
+ {\O(1)\cdot}\int \Big|v+\frac{Rw}{2c}-\frac{Sz}{2c}\Big|\,|R|\, \W^-dx
+{\O(1)\cdot}\int \Big|v+\frac{Rw}{2c}-\frac{Sz}{2c}\Big| |S|\, \W^+dx
\cr\cr
\qquad\qquad 
\ds+ \O(1)\cdot \int |w|\bigl(|RS|+S^2\bigr)\, \W^-\, dx
+ \O(1)\cdot \int |z|\bigl(|RS|+S^2\bigr)\, \W^+\, dx\cr\cr
\ds\qquad\qquad+a(t)\int \left|w_x+\frac{c'}{4c^2}(w-z)S\right|\W^-\, dx
-2 c_0\cdot\int \left|w_x+\frac{c'}{4c^2}(w-z)S\right|\, S^2\,\W^-\, dx
\,.\enda
\eeq

{\bf 5.}
To estimate the time derivative of $I_5$, using (\ref{4.4})  we compute
\bel{4.5}\bega{l}
\ds\left[R\Big(w_x+\frac{c'}{4c^2}(w-z)S\Big)\right]_t
+\left[Rc\Big(w_x+\frac{c'}{4c^2}(w-z)S\Big)
\right]_x\cr\cr
\qquad \ds =~\frac{c'}{2c}R\tilde{s}-\frac{c'}{2c}R\tilde{r}+\frac{c'}{2c}RS
\Big(w_x+\frac{c'}{4c^2}(w-z)S\Big)
\cr\cr\qquad\qquad\ds
+\frac{c'}{2c}R\Big(Sz_x+\frac{c'}{4c^2}(w-z)SR\Big)
-\frac{c'}{2c}R\Big(Rw_x+\frac{c'}{4c^2}(w-z)SR\Big)
\cr\cr\qquad\qquad\ds
-\frac{c''c-(c')^2}{2c^2}R^2\Big(v+\frac{Rw}{2c}-\frac{Sz}{2c}\Big)
+\frac{c''c-2(c')^2}{2c^2}SR\Big(v+\frac{Rw}{2c}-\frac{Sz}{2c}\Big)
\cr\cr\qquad\qquad\ds
+\frac{c''c-(c')^2}{4c^3}R^2S(w-z)-\frac{(c')^2}{8c^3}S^2R(w-z)
\cr\cr\qquad\qquad\ds
+\frac{c'}{4c}(R^2-S^2)\Big(w_x+\frac{c'}{4c^2}(w-z)S\Big)
\cr\cr
\qquad \ds =~\frac{c'}{2c}R\tilde{s}-\frac{c'}{4c}\Big(2R\tilde{r}+R^2w_x
+\frac{c'}{4c^2}(w-z)SR^2\Big)
\cr\cr\qquad\qquad\ds+\frac{c'}{2c}S\Big(Rw_x+\frac{c'}{4c^2}(w-z)RS\Big)
+\frac{c'}{2c}R\Big(Sz_x+\frac{c'}{4c^2}(w-z)SR\Big)\cr\cr\qquad\qquad\ds
-\frac{c''c-(c')^2}{2c^2}R^2\Big(v+\frac{Rw}{2c}-\frac{Sz}{2c}\Big)
+\frac{c''c-2(c')^2}{2c^2}SR\Big(v+\frac{Rw}{2c}-\frac{Sz}{2c}\Big)
\cr\cr\qquad\qquad\ds
+\frac{c''c-(c')^2}{4c^3}R^2S(w-z)-\frac{(c')^2}{8c^3}S^2R(w-z)-\frac{c'}{4c}S^2\Big(w_x+
\frac{c'}{4c^2}(w-z)S\Big).\enda
\eeq
We thus conclude  
\bel{DI5}\bega{l}
\ds\frac{d}{dt}\int  \Big|Rw_x+\frac{c'}{4c^2}(w-z)RS\Big|\,\W^- 
\,dx\cr\cr
\leq~ \ds O(1)\cdot 
\int| R\tilde{s}|\,\W^-\, dx +\O(1)\cdot \int \Big|2R\tilde{r}+R^2w_x
+\frac{c'}{4c^2}(w-z)SR^2\Big|\,\W^-\, dx
\cr\cr
\qquad \ds +\O(1)\cdot\int \Big|Rw_x+\frac{c'}{4c^2}(w-z)RS\Big|\,|S| \W^-dx\cr\cr
\qquad \ds+\O(1)\cdot\int \Big|Sz_x+\frac{c'}{4c^2}(w-z)SR\Big|\, |R|\,\W^-dx 
\cr\cr
\qquad \ds
+ \O(1)\cdot \int \Big|v+\frac{Rw}{2c}-\frac{Sz}{2c}\Big|\,R^2\, \W^-dx
+\O(1)\cdot\int \Big|v+\frac{Rw}{2c}-\frac{Sz}{2c}\Big| |RS|\, \W^+dx
\cr\cr
\qquad 
\ds
+ \O(1)\cdot \int \bigl(|w|+|z|\bigr)(1+R^2)(1+S^2)\, \W^-\, dx
+\O(1)\cdot\int S^2\Big|w_x+\frac{c'}{4c^2}(w-z)S\Big|\, \W^-dx
\cr\cr\qquad\ds
 +a(t)\cdot\int \left|Rw_x+\frac{c'}{4c^2}(w-z)RS\right|\,\W^-\, dx
-2 c_0\cdot\int \left|Rw_x+\frac{c'}{4c^2}(w-z)RS\right|\, S^2\,\W^-\, dx
\enda
\eeq

\v
{\bf 6.} Finally, to 
estimate the time derivative of $I_6$ (change in base measure with density $R^2$), 
we compute
\bel{4.6}
\bega{l}
\ds \left(2R\tilde{r}+ R^2w_x+\frac{c'}{4c^2}(w-z)SR^2\right)_t+\left(c
\Big(2R\tilde{r}+R^2w_x+\frac{c'}{4c^2}(w-z)SR^2\Big)\right)_x
\cr\cr 
\ds\qquad 
=~(R_t-cR_x)\Big(2\tilde{r}+ Rw_x+\frac{c'}{4c^2}(w-z)SR\Big)\cr\cr \ds\qquad 
\quad  +R\,\left[2\Big(\tilde{r}_t-(c\tilde r)_x\Big)+ \Big((Rw_x)_t- (c Rw_x)_x\Big)
+\Big( \frac{c'}{4c^2}(w-z)SR\Big)_t -\Big(\frac{c'}{4c}(w-z)SR)\Big)_x\right]\cr\cr 
\qquad 
 \ds=~\frac{c'}{4c}(R^2-S^2)\Big(2\tilde{r}+ Rw_x+\frac{c'}{4c^2}(w-z)SR\Big)\cr\cr 
 \ds\qquad 
\qquad
+\frac{c''c-{c'}^2}{2c^2}R^3\Big(v+\frac{Rw}{2c}-\frac{Sz}{2c}\Big)-\frac{c''c-2{c'}^2}{2c^2}S^2R\Big(v+\frac{Rw}{2c}-\frac{Sz}{2c}\Big)
\cr\cr \ds\qquad 
\qquad +\frac{c''c-{c'}^2}{4c^3}R^3S(z-w)+\frac{(c')^2}{8c^3}R^2S^2(w-z)
\cr\cr \ds\qquad 
\qquad -\frac{c'}{2c}R\Big(2S\tilde{s}+S^2z_x+\frac{c'}{4c^2}S^2R(w-z)\Big)
+\frac{c'}{c}SR\tilde{r}
\cr\cr \ds\qquad \qquad
+\frac{c'}{2c}R^2\tilde{s}-\frac{c'}{4c}R\Big(2R\tilde{r}+R^2w_x+\frac{c'}{4c^2}(w-z)SR^2\Big)
+\frac{c'}{2c}SR\Big(Rw_x+\frac{c'}{4c^2}(w-z)RS\Big)
\cr\cr \ds\qquad 
\qquad+\frac{c'}{2c}R^2(Sz_x+\frac{c'}{4c^2}(w-z)SR)
\cr\cr \ds\qquad 
\qquad -\frac{c''c-(c')^2}{2c^2}R^3\Big(v+\frac{Rw}{2c}-\frac{Sz}{2c}\Big)
+\frac{c''c-2(c')^2}{2c^2}SR^2\Big(v+\frac{Rw}{2c}-\frac{Sz}{2c}\Big)
\cr\cr \ds\qquad \qquad
+\frac{c''c-(c')^2}{4c^3}R^3S(w-z)-\frac{(c')^2}{8c^3}S^2R^2(w-z)
-\frac{c'}{4c}S^2R\Big(w_x+\frac{c'}{4c^2}(w-z)S\Big)
\cr\cr
\qquad\ds =~\frac{c'}{2c}R^2 \Big(Sz_x+\frac{c'}{4c^2}(w-z)RS\Big)
-\frac{c'}{2c}S^2\Big(Rw_x +\frac{c'}{4c^2} (w-z)SR\Big)\cr\cr\qquad\qquad \ds
+\frac{c''c-2c'^2}{2c^2}\Big(v+\frac{Rw}{2c}-\frac{Sz}{2c}\Big)(R^2S-RS^2)
-\frac{c'}{2c} S^2\tilde{r}+\frac{c'}{2c}R^2\tilde{s}\cr\cr\qquad\qquad \ds
+\frac{c'}{2c}S\left(2R\tilde{r}+R^2w_x+\frac{c'}{4c^2}R^2S(w-z)\right)
-\frac{c'}{2c}R\left(2S\tilde{s}+S^2z_x+\frac{c'}{4c^2}RS^2(w-z)\right).
\enda
\eeq
This yields the estimate
\bel{DI6}
\bega{l}
\ds\frac{d}{dt}\int \Big| 2R\tilde{r}+ R^2w_x+\frac{c'}{4c^2}(w-z)SR^2\Big|\,\W^- \,dx\cr\cr
\qquad\leq~ 
\ds\O(1)\cdot\int R^2 \Big|Sz_x+\frac{c'}{4c^2}(w-z)RS\Big|\,\W^+ dx
\cr\cr
\qquad\qquad \ds
+\O(1)\cdot\int S^2\Big|Rw_x +\frac{c'}{4c^2} (w-z)RS\Big|\,\W^-\, dx
\cr\cr\qquad\qquad \ds
+\O(1)\cdot\int \Big|v+\frac{Rw}{2c}-\frac{Sz}{2c}\Big|\,|R^2S-RS^2|\,\W^-\,dx\cr\cr
\qquad\qquad \ds
+\O(1)\cdot\int S^2|\tilde{r}|\,\W^-\, dx+
\O(1)\cdot\int R^2|\tilde{s}|\,\W^+\, dx\cr\cr\qquad\qquad \ds
+\O(1)\cdot\int |S|\left|2R\tilde{r}+R^2w_x+\frac{c'}{4c^2}R^2S(w-z)\right|\,\W^-\, dx
\cr\cr\qquad\qquad \ds
+\O(1)\cdot\int |R|\left|2S\tilde{s}+S^2z_x+\frac{c'}{4c^2}RS^2(w-z)\right|\,\W^+\, dx
\cr\cr
\qquad\qquad\ds +\int \bigl(a(t)-2 c_0\bigr)
\left| 2R\tilde{r}+ R^2w_x+\frac{c'}{4c^2}(w-z)SR^2\right|\,S^2\,\W^- \,dx\,.
\enda
\eeq
\v

\begin{figure}[htbp]
   \centering
\includegraphics[width=0.7\textwidth]{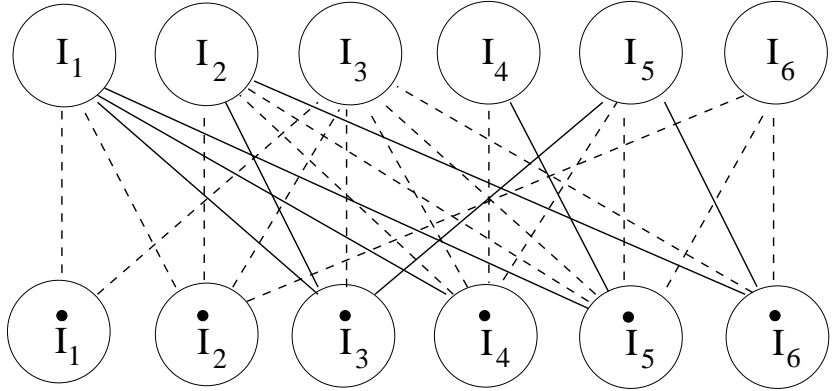}
   \caption{\small A graphical summary of all the a priori estimates.
  If a lower box $\dot I_k$ is connected to an upper box $I_{\ell}$, this 
   means that the integral $I_\ell$ is used in order
   to control the time derivative $\dot I_k ={d\over dt} I_k$. 
   If $\ell\in \F^\sharp_k$, then $\dot I_k$ and $I_\ell$ are connected by a solid line.
     If $\ell\in \F^\flat_k$, then $\dot I_k$ and $I_\ell$ are connected by a dashed line.
 }
   \label{f:wa36}
\end{figure}

{\bf 7.} We keep track of all the above 
estimates by the diagram in Fig.~\ref{f:wa36}.
Recalling (\ref{main}), the weighted norm of a tangent vector can be written as
\bel{norm}\bega{rl}
\Big\|(\tilde{r},w,\tilde{s},z) \Big\|_{(u,R,S)}&=~
\kappa_1 I_1+\kappa_2I_2+\kappa_3I_3+\kappa_4I_4+\kappa_5I_5+ \kappa_6I_6\cr\cr
\ds
\qquad &\ds=~\sum_{k=1}^6 \kappa_k\left(\int J_k^-\,\W^-\, dx+\int J_k^+\,\W^+\, dx
\right),
\enda\eeq
where $J_k^-, J_k^+$ are the various integrands.
According to the estimates (\ref{DI1}), (\ref{DI2}), (\ref{DI3}), (\ref{DI4}), (\ref{DI5}), 
and (\ref{DI6}), the time derivative of each $I_k$ can be estimated as 
\bel{dotik}\bega{rl}
\dot I_k&\ds\leq~
\sum_{\ell\in \F^\flat_k} 
\O(1)\cdot\left(\int J_\ell^-\,\bigl(1+|S|\bigr)\,\W^-\, dx+\int J_\ell^+\, 
\bigl(1+|R|\bigr)\,\W^+\, dx
\right)\cr\cr
&\qquad +\ds\sum_{\ell\in \F^\sharp_k} 
\O(1)\cdot\left(\int J_\ell^-\,(1+ R^2)\,\W^-\, dx+\int J_\ell^+\, (1+S^2)\,\W^+\, dx
\right)\cr\cr
&\qquad +\ds  a(t) I_k - 2c_0\left(\int S^2 J_k^-\,\W^-\, dx+\int R^2J_k^+\,\W^+\, dx
\right).
\enda
\eeq
Here $\F^\flat_k,\F^\sharp_k\subset\{ 1,2,\ldots,6\}$ are suitable sets of indices, illustrated in Fig.~\ref{f:wa36}.
By direct inspection, we see that the set-valued map 
$k\mapsto \F^\sharp_k$
has no cycles.   Indeed, the composition  
$\F^\sharp_k\circ\F^\sharp_k\circ\F^\sharp_k$ yields the empty set.

By choosing a constant $\delta>0$ small enough, we thus obtain a weighted norm
\bel{WN}
\Big\|(\tilde{r},w,\tilde{s},z) \Big\|_{(u,R,S)}~\doteq~
I_1 + \delta I_2 + \delta^3 \,I_3
+\delta I_4 + \delta^2 \,I_5 + \delta^3 \,I_6\eeq
which satisfies the desired inequality (\ref{wn1}).
This completes the proof of Proposition~1.
\endproof
\v

\section{Tangent vectors in transformed coordinates}
\setcounter{equation}{0}
Given any path
$\theta\mapsto u^\theta$, $\theta\in [0,1]$  
of smooth solutions to (\ref{2}),
the analysis in the previous section has
provided an estimate on how its weighted length increases in time.
However, even for smooth
initial data, it is well known that
the quantities $u_t, u_x$ can blow up in finite time
\cite{GHZ}.  When this happens,  a tangent vector 
may no longer exist; even if it does exist,  it is not obvious that
our earlier estimates should remain valid.
Aim of this section is to address these issues.  
Roughly speaking, we claim that
\begi
\item[(i)] Every path of solutions 
$\theta\mapsto u^\theta$
 can be uniformly approximated by a second path 
$\theta\mapsto \tilde u^\theta$
such that, for all but finitely many values of $\theta\in [0,1]$,
 the solution $\tilde u^\theta$ is piecewise smooth, with 
``generic" singularities.
\item[(ii)] If all solutions $u^\theta$ are piecewise smooth, with 
``generic" singularities along  finitely many points and finitely many 
curves in the $t$-$x$ plane, then the tangent vectors are 
still well defined
and their norms can be estimated as before.
\endi

A precise formulation of (i) was recently proved by the authors 
in \cite{BChen}.
The proof is based on the representation of solutions to (\ref{2}) 
in terms of a semilinear system with smooth coefficients
\cite{BZ}, followed by an application of Thom's transversality theorem.
We review here this basic construction, and the characterization
of generic (structurally stable) singularities \cite{BHY}.

To deal with possibly unbounded values of $R,S$ 
in (\ref{0.3}), following \cite{BZ}
it is convenient to introduce a new set of dependent variables:
\bel{abdef}
\alpha~\doteq ~2\arctan R\,,\qquad\qquad \beta~\doteq~ 2\arctan S\,.
\eeq
Using (\ref{0.5}), we obtain the equations
\bel{a0}
\alpha_t-c\,\alpha_x~=~{2\over 1+R^2}(R_t-c\,R_x)
~=~{c'\over 2c}{R^2-S^2\over 1+R^2}\,, 
\eeq
\bel{b0}
\beta_t+c\,\beta_x~=~{2\over 1+S^2}(S_t+c\,S_x)
~=~{c'\over 2c}{S^2-R^2\over 1+S^2}\,. 
\eeq

We now
perform a further change of independent variables.
Consider
the equations for the  backward and forward characteristics:
\beq
\dot x^-~=~-c(u)\,,\qquad\qquad \dot x^+~=~c(u)\,,\label{2.12}
\eeq
where the upper dot denotes a derivative w.r.t.~time. 
The characteristics passing through the point $(t,x)$
will be denoted by
$$
s~\mapsto~ x^-(s,t,x)\,,\qquad\qquad s~\mapsto ~x^+(s,t,x)\,,
$$
respectively.
We shall use a set of coordinates $(X,Y)$ 
on the $t$-$x$ plane such that $X$ is constant along backward characteristics and $Y$ is constant along forward characteristics,
namely
\bel{2.15}\left\{\bega{rl}
X_t- c(u)X_x& =~0\,,\\[4mm]Y_t+ c(u)Y_x &=~0\,.
\enda\right.
\eeq
For example, one can define $X,Y$ to be the intersections
with the $x$-axis,
of the characteristics through the point $(t,x)$, i.e.
\bel{XY1}
X(t,x)~\doteq~ x^-(0,t,x)
\,,\qquad\qquad Y(t,x)~\doteq ~- x^+(0,t,x)\,.
\eeq
More generally, one can consider strictly increasing 
functions 
$x\mapsto \ov X(x)$ and $x\mapsto \ov Y(x)$ and define
\bel{XY2}
X(t,x)~\doteq~\ov X\bigl( x^-(0,t,x)\bigr)
\,,\qquad\qquad Y(t,x)~\doteq ~\ov Y\bigl(- x^+(0,t,x)\bigr)\,.
\eeq

For any smooth function $f$, using (\ref{2.15}) one finds
\beq\left\{
\begin{array}{ccccccr}
f_t+cf_x &=& f_XX_t+f_Y Y_t+cf_X X_x+cf_Y Y_x & = & (X_t+cX_x)f_X
& = & 2cX_x f_X\,,\cr\cr
f_t-cf_x &=& f_XX_t+f_Y Y_t-cf_X X_x-cf_Y Y_x & = & (Y_t-cY_x)f_Y
& = & -2cY_x f_Y\,.
\end{array}\right. \label{2.16} 
\eeq
We now introduce the further variables
\bel{2.17}
p~\doteq~ {1+R^2\over X_x}\,,\qquad\qquad q~\doteq ~{1+S^2\over -Y_x}\,.
\eeq
Notice that the above definitions imply
\bel{2.18}
{1\over X_x}~=~ {p\over 1+R^2}~=~{(1+\cos\alpha)p\over 2}\,,
\qquad\qquad
{-1\over Y_x}~=~{q\over 1+S^2}~=~{(1+\cos\beta)q\over 2}\,.
\eeq

Starting with the nonlinear equation (\ref{0.1}),
using $X,Y$ as independent variables one obtains a semilinear
hyperbolic system with smooth coefficients for the variables $u,\alpha,\beta,p,q$, namely
\bel{u}
\left\{
\begin{array}{ccc}
u_X &= & {\sin \alpha\over 4c} \, p\,,\\[4mm]
 u_Y &= & {\sin \beta\over 4c} \, q\,,
\end{array}\right.
\eeq
\bel{ab}
\left\{
\begin{array}{ccc}
\alpha_Y&=&{c'\over 8c^2}\,( \cos \beta - \cos \alpha)\,q\,,\\[4mm]
\beta_X&=&{c'\over 8c^2}\,( \cos \alpha - \cos \beta)\,p\,,
\end{array}\right. 
\eeq
\bel{pq}
\left\{
\begin{array}{ccc}
p_Y &= & {c'\over 8c^2}\, \big(\sin \beta-\sin \alpha\big)\,pq\,,\\[4mm]
q_X &= & {c'\over 8c^2}\, \big(\sin \alpha-\sin \beta\big)\,pq\,.
\end{array}\right.
\eeq
The map $(X,Y)\mapsto (t,x)$ can  be constructed as follows.
Setting $f=x$, then $f=t$  in the two equations at
(\ref{2.16}), we find
$$
\left\{
\begin{array}{rcr}
c&=& 2cX_x\,x_X\,, \\[3mm]
-c&=& -2cY_x\,x_Y\,,
\end{array}\right.\qquad\qquad
\left\{
\begin{array}{rcr}
1&= &2c\,X_x\,t_X\,,\\[3mm]
1&=& -2c\,Y_x\,t_Y\,,
\end{array}\right.
$$
respectively.
Therefore, using (\ref{2.18}) we obtain
\bel{x}
\left\{
\begin{array}{ccccr}
x_X&~=~&{1\over 2X_x}&~=~&{(1+\cos \alpha)\,p\over 4}\,,\\[4mm]
x_Y&~=~&{1\over 2Y_x}&~=~&-{(1+\cos \beta)\,q\over 4}\,,
\end{array}\right. 
\eeq
\bel{t}
\left\{
\begin{array}{ccccr}
 t_X&=&{1\over 2cX_x}&=&{(1+\cos \alpha)\,p\over 4c}\,,\\[4mm]
t_Y&=&{1\over -2cY_x}&=&{(1+\cos \beta)\,q\over 4c}\,.
\end{array}\right. \eeq
\v
Given the initial data (\ref{0.2}), one particular way to assign the corresponding  boundary 
data for   (\ref{u})-(\ref{t}) is as follows.
In the $X$-$Y$ plane, consider the line 
\bel{g0def}\gamma_0~=~\{X+Y=0\}~\subset~\R^2\eeq
parameterized as
$x~\mapsto ~(\ov X(x), \,\ov Y(x))~\doteq~(x,\, -x)$.
Along $\gamma_0$ we can assign the boundary data
$(\ov u,\ov \alpha, \ov \beta, \ov p, \ov q)$ by setting
 \beq{\ov u }~ =~u_0(x)\,,\qquad\qquad
\left\{
\begin{array}{rcl}
{\ov \alpha } &= & 2\arctan R(0,x)\,,\\ {\ov \beta} &= & 2\arctan S(0,x)\,,
\end{array}\right.\qquad \qquad
\left\{
\begin{array}{rcl}
{\ov p} &\equiv & 1+ R^2(0,x)\,,\\
{\ov q} &\equiv & 1+S^2(0,x)\,,
\end{array}\right.  \label{2.28}
\eeq
at each point $(x,-x)\in\gamma_0$.
We recall that, at time $t=0$, by (\ref{0.2}) one has
$$\bega{lr}R(0,x) &=~(u_t + c(u) u_x)(0,x) ~=~
u_1(x) + c(u_0(x)) u_{0,x}(x),\\[4mm]
S(0,x) &=~(u_t - c(u) u_x)(0,x) ~=~u_1(x) - c(u_0(x)) u_{0,x}(x).\enda$$
\v
{\bf Remark 3.}  The above construction (\ref{g0def})-(\ref{2.28})
is by no means the unique way to prescribe initial values.
One should be aware that  many distinct solutions to the system
(\ref{u})--(\ref{t}) can yield the same solution $u=u(t,x)$ of (\ref{0.1})-(\ref{0.2}).
Indeed, let $(u,\alpha,\beta ,p,q,x,t)(X,Y)$ be one particular solution. 
Let $\phi,\psi:\R\mapsto\R$ be two $\C^2$ bijections, 
with $\phi'>0$ and $\psi'>0$.  Introduce the new independent and dependent variables
$(\Tilde X,\Tilde Y)$ and $(\tilde u,\tilde \alpha,\tilde \beta ,
\tilde p,\tilde q,\tilde x,\tilde t)$
by setting
\bel{TXY}
X~=~\phi(\Tilde X)\,,\qquad\qquad Y~=~\psi(\Tilde Y),\eeq
\bel{TUWZ}
(\tilde u,\tilde \alpha,\tilde \beta ,\tilde x, \tilde t)(\Tilde X,\Tilde Y)~=~(u,\alpha,\beta ,x,t)(X,Y),\eeq
\bel{TPQ}\left\{\bega{rl}
\tilde p(\Tilde X,\Tilde Y)&=~p(X,Y)\cdot \phi'(\Tilde X),\\[4mm]
\tilde q(\Tilde X,\Tilde Y)&=~q(X,Y)\cdot \psi'(\Tilde Y).\enda\right.\eeq
Then, as functions of $(\Tilde X, \Tilde Y)$,  the
variables $(\tilde u,\tilde \alpha,\tilde \beta ,\tilde p,\tilde q,\tilde x,\tilde t)$
provide another solution of the same system (\ref{u})--(\ref{t}). 
Moreover, by (\ref{TUWZ}) the set
\bel{graph2}\Big\{ \Bigl(\tilde t(\Tilde X,\Tilde Y), \, \tilde x(\Tilde X,\Tilde Y), 
\,\tilde u(\Tilde X,\Tilde Y)\Bigr)\,;~~
(\Tilde X,\Tilde Y)\in\R^2\Big\}\eeq
coincides with the set in (\ref{graph}). Hence it is the graph of the same solution 
$u$ of (\ref{0.1}).  
One can regard the variable transformation (\ref{TXY})
simply as a relabeling of forward and backward characteristics, in the solution $u$. In 
connection with first order wave  equations,  relabeling symmetries
have been studied in
\cite{BHR, GHR2}.
\v
{\bf Remark 4.}  The  system (\ref{u})--(\ref{t}) is clearly 
invariant w.r.t.~the addition
of an integer multiple of $2\pi$ to the variables $\alpha,\beta $.
Taking advantage of this property, in the following we 
shall regard $\alpha,\beta $ as points in the quotient manifold
${\mathbb T}\doteq \R/2\pi {\mathbb Z}$.  
As a consequence, we have the implications
\bel{impl}\bega{rl}
\alpha&\not=~\pi\qquad\implies\qquad \cos \alpha~>~-1\,,\\[3mm]
\beta &\not=~\pi\qquad\implies\qquad \cos \beta ~>~-1\,.\enda\eeq
\v
{\bf Remark 5.}
Since the semilinear system (\ref{u})--(\ref{t}) has smooth coefficients,
for smooth initial data all components of the solution remain 
smooth on the entire $X$-$Y$
plane.   As proved in \cite{BZ}, 
the quadratic terms in (\ref{pq}) (containing the product
$pq$)  account
for transversal wave interactions and 
do not produce finite time blowup of the
variables $p,q$.   Moreover, if the values of $p,q$ are uniformly 
positive and bounded on the  line  $\gamma_0$, then 
they remain uniformly positive and bounded
on compact sets of the $X$-$Y$ plane.  Throughout this paper, we always
consider solutions of (\ref{u})--(\ref{t}) where $p,q>0$.
\v
The main results in \cite{BCZ, BZ} can be summarized as
\v
{\bf Theorem~3.} {\it  Let $c=c(u)$ be a smooth, uniformly positive function.
Let  $(t,x,u,\alpha,\beta,p,q)(X,Y)$ be a smooth solution of 
the semilinear system
(\ref{u})--(\ref{t}) with boundary data as in (\ref{2.28}). Then 
the function $u=u(t,x)$ whose graph is
\bel{graph}
\hbox{Graph}(u)~=~\Big\{\bigl(t(X,Y), \, x(X,Y),\, u(X,Y)
\bigr)\,;~~(X,Y)\in \R^2\Big\}\eeq
provides the unique conservative solution to the Cauchy problem
(\ref{0.1})-(\ref{0.2}).}
\v

\begin{figure}[htbp]
   \centering
\includegraphics[width=0.6\textwidth]{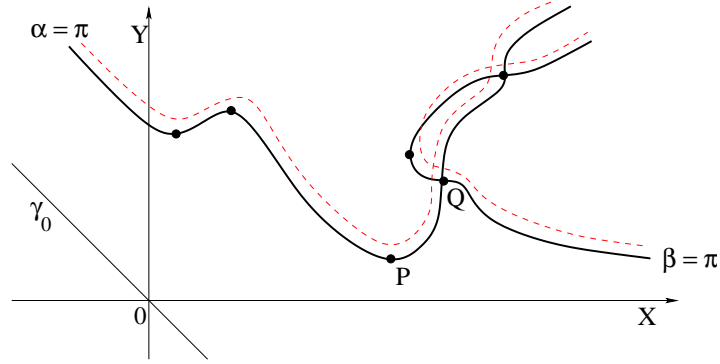}
   \caption{\small The level sets $\{\alpha=\pi\}$ and $\{\beta=\pi\}$
   in a solution with generic singularities.  In the $X$-$Y$ plane
   these
   are smooth curves which are structurally stable 
   w.r.t.~small $\C^2$ perturbations.
 }
   \label{f:wa83}
\end{figure}

Throughout the following, we shall be interested not in a single solution but in a continuous path of solutions $\theta\mapsto u^\theta$, $\theta\in [0,1]$.   
We introduce
suitable regularity conditions, allowing us to compute the ``length" of this path
by integrating a suitable norm of its tangent vector $\|du^\theta(t,\cdot)/ d\theta\|$.
\v
{\bf Definition 1.}  {\it  We say that a solution $u=u(t,x)$ of 
(\ref{0.1}) has {\bf generic singularities} for $t\in [0,T]$
if it admits a representation of the form (\ref{graph}), where
(i) the functions $(u,\alpha,\beta,p,q,x,t)(X,Y)$ are $\C^\infty$, 
and  (ii)  on the domain where  $t(X,Y)\in [0,T]$  the following generic conditions
hold:
\begi
\item[(G1)] \qquad $\alpha=\pi$,~ $\alpha_X = 0\qquad\implies\qquad
\alpha_Y\not= 0$,~ $\alpha_{XX}\not= 0$,

\item[(G2)] \qquad $\beta=\pi$, ~$\beta_Y = 0\qquad\implies\qquad
\beta_X\not= 0$, ~$\beta_{YY}\not= 0$,

\item[(G3)] \qquad$\alpha=\pi$, ~$\beta=\pi, \qquad\implies\qquad
\alpha_X\not= 0$,~ $\beta_Y\not= 0$.
\endi
}

\begin{figure}[htbp]
   \centering
\includegraphics[width=0.5\textwidth]{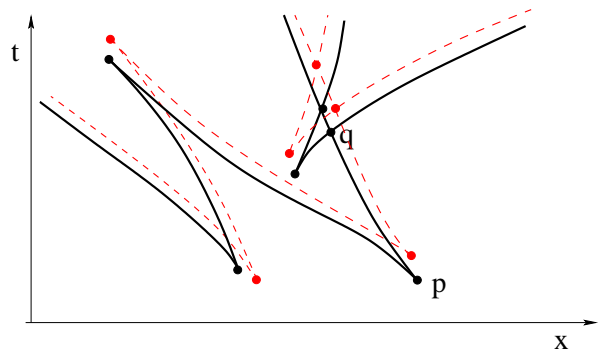}
   \caption{\small The set of singular points  (where $|u_x|\to+\infty$)
   in a solution $u(t,x)$. These are the images of the sets
   $\{\alpha=\pi\}$ and $\{\beta=\pi\}$ in Fig.~\ref{f:wa83}. 
   By structural stability, 
   every small perturbation will yield anther solution with the same
   type of singularities.
 }
   \label{f:wa82}
\end{figure}

\begin{figure}[htbp]
\centering
\includegraphics[width=1.0\textwidth]{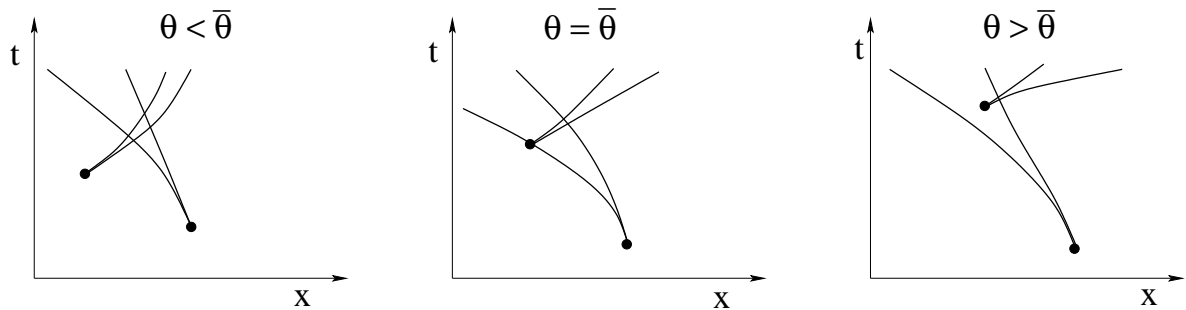}
\caption{ \small 
Here the solution $u^\theta$ has generic (i.e., structurally stable) singularities for $\theta<\bar \theta$ and for $\theta>\bar\theta$.
However, when the parameter $\theta$ crosses the critical 
value $\bar \theta$, the topology of the singular set changes.}
\label{f:w27}
\end{figure}

Some words of explanation are in order.  
Even if the solution $(X,Y)\mapsto (x,t,u,\alpha,\beta ,p,q)(X,Y)$ of the semilinear system
(\ref{u})--(\ref{t})
remains smooth on the entire $X$-$Y$ plane, 
the function
$u=u(t,x)$ in (\ref{graph}) can have singularities because the 
coordinate change $\Lambda:(X,Y)\mapsto (x,t)$ is not smoothly invertible. Indeed,
by (\ref{t})-(\ref{x}), the Jacobian matrix is computed by
\bel{J}D\Lambda~=~
\left(\bega{cc} x_X & x_Y\cr
t_X& t_Y\enda\right)~=~
\left(\bega{ccc} {(1+\cos \alpha)\,p\over 4} && -{(1+\cos \beta )\,q\over 4}\\[3mm]
{(1+\cos \alpha)\,p\over 4c(u)}&& {(1+\cos \beta )\,q\over 4c(u)}\enda\right)
\eeq
We recall that $p,q$ remain uniformly positive and uniformly bounded
on compact subsets of the $X$-$Y$ plane.
By Remark~3, at a point $(X_0,Y_0)$ where  
$\alpha\not=\pi$ and $\beta \not=\pi$, this  matrix
is  invertible, having a strictly positive determinant.   
The function $u=u(x,t)$ considered at 
(\ref{graph}) is thus smooth on a neighborhood of the point 
$$( t_0, x_0) ~=~\bigl( t(X_0, Y_0)\,,~x(X_0,Y_0)\Bigr).$$
To study the set of points in the $x$-$t$ plane where $u$ is singular,
we thus need to look at points where either $w=\pi$ or $\beta =\pi$.
The generic conditions (G1)--(G2) guarantee that these level sets
are smooth curves in the $X$-$Y$ plane.
Condition (G3) implies that the level sets where $\{\alpha=\pi\}$
and $\{\beta=\pi\}$ intersect transversally because $\alpha_Y=\beta_X=0$
when $\alpha=\beta=0$.
As observed in \cite{BChen}, the conditions (G1)--(G3) 
are invariant w.r.t.~smooth variable transformations 
$(X,Y)\leftrightarrow (\Tilde X, \Tilde Y)$.
We also remark that, if a solution $U=(u,\alpha,\beta,p,q)$
of (\ref{u})--(\ref{pq})  satisfies the generic conditions 
(G1)--(G3),   then 
by the implicit function theorem  the same is true for every
perturbed solution 
$\Tilde U=(\tilde u,\tilde \alpha,\tilde \beta,\tilde p,\tilde q)$  sufficiently close to $U$. In other words, generic singularities are
{\em structurally stable}.  An example of structurally unstable solution,
corresponding to a change of topology in the singular set,
is shown in Fig.~\ref{f:w27}. 
\v
{\bf Definition 2.}  {\it 
 We say that a path of initial data 
 $\gamma:\theta\mapsto 
 (u_0^\theta, 
 u_1^\theta)$, $\theta \in [0,1]$ is a {\bf piecewise regular path}
if the following conditions are satisfied.
\begi
\item[(i)] There exists a continuous 
map $(X,Y,\theta)\mapsto (u,\alpha,\beta,p,q,x,t)$ such that, for each
$\theta\in [0,1]$ the semilinear system (\ref{u})--(\ref{t})
is satisfied. Moreover, the function $u^\theta(x,t)$ whose graph is
\bel{gh2}
\hbox{Graph}(u^\theta)~=~\Big\{(t,x,u)(X,Y,\theta);~~(X,Y)\in \R^2
\Big\}\eeq
provides the conservative solution of (\ref{2}) with initial data
$$u^\theta(0,x)~=~u^\theta_0(x)\,,
\qquad\qquad u^\theta_t(0,x)~=~u^\theta_1(x)\,.$$
\item[(ii)] There exist finitely many values
$0=\theta_0<\theta_1<\cdots <\theta_N = 1$ such that the following holds. For $\theta\in \,]\theta_{i-1},\theta_i[$, the map
$(X,Y,\theta)\mapsto (u,\alpha,\beta,p,q,x,t)$ is $\C^\infty$. 
Moreover, 
the solution $u^\theta=u^\theta(t,x)$ has generic 
singularities at time 
$t=0$.
\endi 

In addition, if for all  $\theta\in [0,1]\setminus\{\theta_1,\ldots,\theta_N\}$, 
the  solution 
$u^\theta$ has generic singularities for $t\in [0,T]$, 
then we say that 
the path of solutions $\gamma:\theta\mapsto u^\theta$ is 
{\bf piecewise regular for} $t\in [0,T]$.
}
\v
{\bf Remark 6.}  According to Remark~3, there are infinitely many
parameterizations of the variables $(X,Y)$ that yield the same 
solution $u=u(t,x)$. However, as shown in \cite{BChen}, 
the property of having generic singularities is 
independent of the particular representation
used in (\ref{gh2}).  
\v
{\bf Remark 7.} The above definition has a simple motivation.
If $\gamma$ is a {\em piecewise regular path}, then we can 
compute its length as an integral of the norm of a tangent 
vector.    In addition, if $\gamma$ is {\em piecewise regular for} 
$t\in [0,T]$,
then  the length of the path of solutions $\gamma^t:\theta\mapsto 
(u^\theta(t,\cdot), u_t^\theta(t,\cdot))$ is well defined
not only at $t=0$ but for all $t\in [0,T]$.  See Definition~3
in Section~6 for details.
\v
{\bf Remark 8.}  In Definition~2,
the finitely many values of $\theta$ where $u^\theta$ does not have structurally stable singularities 
correspond to  bifurcation values.
As  $ \theta$ crosses one of these values, 
the topological structure of the singular set (where
$u^\theta_x\to\pm \infty$) usually changes, as shown
in Fig.~\ref{f:w27}.  

Following \cite{BChen}, on the wave speed $c$ we assume
\begi
\item[{\bf (A)}] The map
 $c:\R\mapsto \R_+$ is
 smooth and  uniformly positive.   The quotient   $c'(u)/c(u)$ is uniformly 
 bounded. Moreover,
 the following generic condition is satisfied:
 \bel{morse}  c'(u)~=~0\qquad\implies\qquad c''(u)~\not=~ 0.\eeq 
 \endi
 Notice that, by (\ref{morse}), the derivative $c'(u)$ 
 vanishes only at isolated points.  The following  result, proved in \cite{BChen},
shows that the set of piecewise regular paths is dense.
\v
{\bf Theorem 4.} 
{\it Let the wave speed $c(u)$ satisfy the assumptions {\bf (A)}
and let $T>0$ be given.
Let $\theta\mapsto (t^\theta, x^\theta, u^\theta, \alpha^\theta, 
\beta^\theta, p^\theta,q^\theta)$, $\theta\in [0,1]$,   
be a smooth path of 
solutions to (\ref{u})--(\ref{t}).  Then there exists
a sequence of  paths of solutions 
$\theta\mapsto (t^\theta_n, x^\theta_n, u^\theta_n, \alpha^\theta_n, 
\beta^\theta_n, p^\theta_n,q^\theta_n)$ with the following properties.
\begi
\item[(i)] For each $n\geq 1$, the  path of corresponding solutions of (\ref{0.1})
$\theta\mapsto u^\theta_n$ is  regular  for $t\in [0,T]$, 
according to Definition~2.
\item[(ii)] For any bounded domain $\Omega$ in the $X$-$Y$ plane,
 as $n\to\infty$ the 
functions $(t^\theta_n, x^\theta_n, u^\theta_n,
 \alpha^\theta_n, 
\beta^\theta_n, $ $ p^\theta_n,q^\theta_n)$ converge to 
$(t^\theta, x^\theta, u^\theta, \alpha^\theta, 
\beta^\theta, p^\theta,q^\theta)$  uniformly in $\C^k([0,1]\times 
\Omega)$, for every $k\geq 1$.
\endi}
\v
Thanks to this density result, to construct a Lipschitz metric
it now remains to 
show that the weighted length of a regular path satisfies
the same  estimates as the smooth paths considered 
in the previous section.
Toward this goal, we first derive an expression for the norm 
of a tangent vector as a line integral in the $X$-$Y$ coordinates.

Consider a reference solution $u$ (\ref{0.1}) and a family of perturbed solutions $u^\ve$, $\ve\in [0, \ve_0]$.
We assume that, in the $X$-$Y$ coordinates, these define a 
smooth family of solutions of (\ref{u})--(\ref{t}), say
$(t^\ve, x^\ve, u^\ve, \alpha^\ve, \beta^\ve, p^\ve,q^\ve)$.
For each $\ve$, the curves where $X$ =constant and $Y$ = constant 
correspond respectively to backward and forward characteristics of the solutions $u^\ve(t,x)$.
We remark that, at time $t=0$, we have considerable freedom in 
choosing these parameterizations.  We can take advantage of this in the following way.
Let $w,z$ be the shifts in (\ref{wzt}).   
At time $t=0$ we choose the parameterizations
according to
\bel{repar}
X^\ve(0,~x+\ve \, w(0,x))~=~x\,,\qquad\qquad  
Y^\ve(0,~x+\ve \, z(0,x))~=~-x\,.\eeq
Consider  the curve in $X$-$Y$ space 
\bel{Gtau}\Gamma_\tau~=~\{(X,Y)\,,~t(X,Y)=\tau\}~
=~\{ (X, Y(\tau,X))\,;~~X\in\R\}~=~\{ (X(\tau, Y),\,Y)\,;~~Y\in\R\}\,,\eeq 
and 
denote by 
\bel{Gtep}\Gamma^\ve_\tau~=~\{(X,Y)\,,~t^\ve(X,Y)=\tau\}~
=~\{ (X, Y^\ve(\tau,X))\,;~~X\in\R\}~=~\{ (X^\ve(\tau, Y),\,Y)\,;~~Y\in\R\}\eeq
the perturbed curve. 
We can write the  perturbed solutions as
\bel{fop}
(t^\ve, x^\ve,u^\ve, \alpha^\ve, \beta^\ve, p^\ve,q^\ve)~=~
(t,x,u, \alpha, \beta, p,q)+\ve ({\mathcal T}, {\mathcal X}, U,A,B,P,Q)+o(\ve)\eeq
Since the system (\ref{t})--(\ref{u})
has smooth coefficients, the first order perturbations satisfy a linearized system
and are well defined for all $(X,Y)\in \R^2$.
We observe that the quantities $v,\tilde r,\tilde s, w,z$ appearing in (\ref{main})
can be expressed in terms of the first order perturbations 
$({\mathcal T}, {\mathcal X}, U,A,B,P,Q)$.
Indeed,
$$(1+R^2)\, dx ~=~p\, dX\,,\qquad\qquad (1+S^2)\, dx~=~-q\, dY$$

Notice that, by definition,
$$t^\ve(X, \, Y^\ve(\tau,X))~=~t^\ve(X^\ve(\tau,Y), Y)~=~\tau.$$
Hence by the implicit function theorem, at $\ve=0$:
$$
{\partial X^\ve\over\partial \ve} ~=~-{\partial t^\ve\over\partial \ve}\cdot 
\left({\partial t\over\partial X}\right)^{-1}~=~-\T\frac{4c}{(1+\cos \alpha )p}$$
and
$${\partial Y^\ve\over\partial \ve} ~=~-{\partial t^\ve\over\partial \ve}\cdot 
\left({\partial t\over\partial Y}\right)^{-1}~=~-\T\frac{4c}{(1+\cos \beta )q}\,.$$

\paragraph{1.}
The shift in $x$ is computed by
$$
\left.
\begin{array}{rcl}
w&=&\lim_{\ve\to 0}~{ x^\ve(X,Y^\ve(\tau,X))- x(X,Y(\tau,X))\over\ve}\\[2mm]
&=&{\mathcal X}(X,Y(\tau,X))
+x_Y\cdot {\partial Y^\ve\over\partial\ve}\Big|_{\ve=0}~=~({\mathcal X}+c{\mathcal T})(X,Y(\tau,X)).
\end{array}\right.
$$

In a similar way,
$$
\left.
\begin{array}{rcl}
z&=&\lim_{\ve\to 0}~{ x^\ve(X^\ve(\tau,Y),Y)- x(X(\tau,Y),Y)\over\ve}\\[2mm]
&=&{\mathcal X}(X(\tau,Y),Y)
+x_X\cdot {\partial X^\ve\over\partial\ve}\Big|_{\ve=0}~=~({\mathcal X}-c{\mathcal T})(X(\tau,Y),Y),
\end{array}\right.
$$

\paragraph{2.}
We now derive an expression for $\tilde r, \tilde s$.  One has
\bel{trs}
r+w R_x~=~{d\over d\ve} \tan{\alpha^\ve(X,Y^\ve(\tau,X))\over 2}\Big|_{\ve=0}
~=~\frac{1}{2}\Big( A-\T\frac{4c}{(1+\cos\beta)q}\alpha_Y \Big)\sec^2\frac{\alpha}{2}
\eeq
and
\bel{trs2}
s+zS_x~=~{d\over d\ve} \tan{\beta^\ve(X^\ve(\tau, Y), Y)\over 2}\Big|_{\ve=0}
~=~\frac{1}{2}\Big( B-\T\frac{4c}{(1+\cos\alpha)p}\beta_X \Big)\sec^2\frac{\beta}{2}\,.
\eeq
By (\ref{rstt}) it thus follows
\bel{trs3}
\tilde{r}
~=~\frac{1}{2}\Big( A-\T\frac{4c}{(1+\cos\beta)q}\alpha_Y \Big)\sec^2\frac{\alpha}{2}
-\frac{c'}{4c}\T\tan^2 \frac{\beta}{2} 
\eeq
and
\bel{trs4}
\tilde{s}
~=~\frac{1}{2}\Big( B-\T\frac{4c}{(1+\cos\alpha)p}\beta_X \Big)\sec^2\frac{\beta}{2}
-\frac{c'}{4c}\T\tan^2 \frac{\alpha}{2} \,.
\eeq

\paragraph{3.}
By (\ref{u}) one has
\[
v+u_x w~=~{d\over d\ve} u^\ve(X, Y^\ve(\tau, X))\Big|_{\ve=0}
~=~U-u_Y \T\frac{4c}{(1+\cos \beta)q}~=~U-\T\tan \frac{\alpha}{2}\,.
\]
Therefore
\beq\label{vUT}
v+\frac{Rw}{2c}-\frac{Sz}{2c}
~=~U-(\tan \frac{\alpha}{2}+\tan \frac{\beta}{2})\cdot 
\T\,.
\eeq

\paragraph{4.}
We now calculate the terms $I_4$ -- $I_6$ in (\ref{main}).

The change in base measure with density $1+R^2$ is given by
\bel{c_1R2}
\lim_{\ve\to 0}~{ p^\ve(X,Y^\ve(\tau,X))- p(X,Y(\tau,X))\over\ve}~=~P(X,Y)
+p_Y\cdot {\partial Y^\ve\over\partial\ve}\Big|_{\ve=0}~=~P- \T\frac{4c}{(1+\cos \beta )q}p_Y\,.
\eeq
The change in base measure with density $1+S^2$ is given by
\bel{c_1S2}
\lim_{\ve\to 0}~{ q^\ve(X^\ve(\tau,Y),Y)- q(X(\tau,Y),Y)\over\ve}~=~Q(X,Y)
+q_X\cdot {\partial X^\ve\over\partial\ve}\Big|_{\ve=0}~=~Q- \T\frac{4c}{(1+\cos \alpha)p}q_X\,.
\eeq

The change in base measure with density $R^2$  (the integrand in $I_6$) is estimated by
\bel{c_R2}
\left.
\begin{array}{l}\ds
{d\over d\ve}\Big( p^\ve \sin^2\frac{\alpha^\ve}
{2}(X,Y^\ve(\tau,X))\Big)\Big|_{\ve=0}\\[4mm]\ds\qquad 
=~\Big(P- \T\frac{4c}{(1+\cos \beta )q}p_Y\Big)\sin^2\frac{\alpha}{2}
~+~\frac{p\sin \alpha}{2}\Big( A-\T\frac{4c}{(1+\cos\beta)q}\alpha_Y \Big)
\,.
\end{array}
\right.
\eeq
%

The difference between (\ref{c_1R2}) and (\ref{c_R2}) shows that the change in base measure with density~1  (the integrand in $I_4$) 
is computed by
\bel{c_1}
\Big(P- \T\frac{4c}{(1+\cos \beta )q}p_Y\Big)\cos^2\frac{\alpha}{2}
~-~\frac{p\sin \alpha}{2}\Big( A-\T\frac{4c}{(1+\cos\beta)q}\alpha_Y \Big)
\,.
\eeq

\v
Combining the previous computations, 
the weighted norm of a tangent vector (\ref{main}) can be written  
as a line integral over the line $\Gamma_\tau$ 
defined at (\ref{Gtau}): 
\bel{N2}
\Big\|(\tilde{r},w,\tilde{s},z) \Big\|_{(u,R,S)}
~=~\sum_{\ell=1}^6 \kappa_\ell \cdot\int_{\Gamma_\tau}
\Big\{|J_\ell| \, \W^-\,dX + |H_\ell|\, \W^+\, dY\Big\}\,,\eeq
where
$$\bega{rcl}
J_1& =& (\X-c\T) p\\[4mm]
J_2&=&\frac{1}{2}\Big( Ap-\T\frac{4cp}{(1+\cos\beta)q}\alpha_Y \Big)
-\frac{c'}{4c}p\T\tan^2 \frac{\beta}{2} \cos^2\frac{\alpha}{2}\\[4mm]
J_3&=&\Big(U-(\tan \frac{\alpha}{2}+\tan \frac{\beta}{2})\cdot \T\Big)p\\[4mm]
J_4&=&P\cos^2\frac{\alpha}{2}- \T\frac{2c}{q}p_Y\frac{\cos^2\frac{\alpha}{2}}{\cos^2\frac{\beta}{2}}
-\frac{p\sin \alpha}{2} A+\frac{2cp}{q}\,\T \alpha_Y\frac{\sin\alpha}{1+\cos \beta}
+\frac{c'}{2c}p\T\tan\frac{\beta}{2}\cos^2\frac{\alpha}{2}\\[4mm]
J_5&=&J_4\cdot\tan\frac{\alpha}{2}\\[4mm]
&=&\frac{1}{2}P\sin\alpha- \T\frac{2c}{q}p_Y\frac{\sin\frac{\alpha}{2}\cos\frac{\alpha}{2}}{\cos^2\frac{\beta}{2}}
-pA\sin^2 \frac{\alpha}{2} +\frac{2cp}{q}\,\T \alpha_Y\frac{\sin^2\frac{\alpha}{2}}{\cos^2\frac{\beta}{2}} 
+\frac{c'}{4c}p\T\tan\frac{\beta}{2}\sin{\alpha} 
\\[4mm]
J_6&=&\Big(P- \T\frac{4c}{(1+\cos \beta )q}p_Y\Big)\sin^2\frac{\alpha}{2}
+\frac{p\sin \alpha}{2} A-\frac{2cp}{q}\,\T \alpha_Y\frac{\sin\alpha}{1+\cos\beta}\\[4mm]
&&\qquad+\frac{c'}{2c}(\sin^2\frac{\alpha}{2}\tan\frac{\beta}{2}-\tan^2\frac{\beta}{2}\sin\frac{\alpha}{2}\cos\frac{\alpha}{2})\T p\,. 
\enda
$$
Using (\ref{pq}) and (\ref{ab}), the above expression can be simplified as 
\bel{J16}
\left\{\bega{rcl}
J_1& =& (\X-c\T) p\\[4mm]
J_2&=&\frac{1}{2} Ap-\frac{c'}{4c}p\T\sin^2\frac{\alpha}{2}\\[4mm]
J_3&=&\Big(U-(\tan \frac{\alpha}{2}+\tan \frac{\beta}{2})\cdot \T\Big)p\\[4mm]
J_4&=&P\cos^2\frac{\alpha}{2}-\frac{p\sin \alpha}{2} A+\frac{c'}{4c}\T p \sin \alpha\\[4mm]
J_5&=&
\frac{1}{2}P\sin\alpha-pA\sin^2 \frac{\alpha}{2}  
+\frac{c'}{2c}\T p \sin^2\frac{ \alpha}{2}
\\[4mm]
J_6&=&P\sin^2\frac{\alpha}{2}
+\frac{p\sin \alpha}{2} A\,.
\enda\right.
\eeq
In a similar way, we obtain
\bel{H16}
\left\{\bega{rcl}
H_1&=&(\X+c\T) q\\[4mm]
H_2&=&\frac{1}{2}Bq -\frac{c'}{4c}q\T\sin^2\frac{\beta}{2}\\[4mm]
H_3&=&\Big(U-(\tan \frac{\alpha}{2}+\tan \frac{\beta}{2})\cdot 
\T\Big)q\\[4mm]
H_4&=&Q\cos^2\frac{\beta}{2}-\frac{q\sin \beta}{2} 
B+\frac{c'}{4c}\T q \sin \beta\\[4mm]
H_5&=&
\frac{1}{2}Q\sin\beta
-qB\sin^2 \frac{\beta}{2}+\frac{c'}{2c}\T q 
\sin^2\frac{ \beta}{2}\\[4mm]
H_6&=& Q\sin^2\frac{\beta}{2}
+\frac{q\sin \beta}{2} B \enda\right.
\eeq
It is clear that the integrands $J_\ell$,  
$H_\ell$ are smooth, for $\ell=1,2,4,5,6$.
We claim that the integrands  $J_3$ and $H_3$ are continuous as well.
 Indeed, using
 \eqref{vUT} we obtain
$$\bega{l}
U-(\tan \frac{\alpha}{2}+\tan \frac{\beta}{2})\cdot \T\\[4mm]
\qquad =~2c v+Rw-Sz\\[4mm]
\qquad =~\int \left(\frac{c'}{c}(R-S)v+2c v_x +wR_x+R w_x- z S_x- S z_x\right)~dx\\[4mm]
\qquad =~\int\Big(r-s+wR_x+R w_x- z S_x- S z_x\Big)\,dx\\[4mm]
\qquad =~\int\Big(r+R w_x-\frac{c'}{8c^2}(w-z)S^2\Big)\,dx -\int\Big(s+zS_x-\frac{c'}{8c^2}(w-z)R^2) \Big)\,dx\\[4mm]
\quad\qquad +\int\Big(wR_x +\frac{c'}{4c^2}R^2S(w-z)\Big)\,dx-\int
\Big(S z_x+\frac{c'}{4c^2}R^2S(w-z)\Big)\,dx\\[4mm]
\qquad\quad +\int\Big(\frac{c'}{8c^2}(w-z)(S^2-R^2)\Big)\,dx\,.
\enda
$$
The three terms on the right hand side 
correspond to the integrands in  $I_2$, $I_4$ and $I_1$,
respectively.  Hence they are continuous.

\section{Length of piecewise regular paths}
\setcounter{equation}{0}

Let $\gamma:\theta\mapsto (u_0^\theta, u_1^\theta)$  be a piecewise 
regular path of initial data. According to Definition~2.
there exists a smooth path of solutions of (\ref{u})--(\ref{t}), say
$\theta\mapsto (x^\theta, t^\theta, u^\theta, \alpha^\theta, 
\beta^\theta, p^\theta,q^\theta)(X,Y)$,
such that (\ref{gh2}) holds for every $\theta\in [0,1]$.
At time $t=0$, an upper bound on the length of this path
can be computed as follows. For each $\theta\in [0,1]$,
consider the curve in the $X$-$Y$ plane
$$\Gamma^\theta_0~\doteq~\Big\{ (X,Y)\,;~~t^\theta(X,Y) =0\Big\}.
$$
The norm of the tangent vector is then determined by (\ref{N2}).
Integrating w.r.t.~$\theta$ we obtain 
\bel{ulg}
\int_0^1\left(\sum_{\ell=1}^6 \kappa_\ell 
\cdot\int_{\Gamma_0^\theta}\Big\{ 
|J_\ell^\theta| \, \W^- dX~+ |H_\ell^\theta|\,\W^+\, dY\Big\}\right)
\, d\theta\,.\eeq

We recall  that  there exist
infinitely many paths of solutions of (\ref{u})--(\ref{t})
which yields the same path of solutions to (\ref{0.1}).
Indeed, as shown in Remark~3, at time $t=0$ for each $\theta$ 
one can choose smooth,
increasing  functions $\phi^\theta, \psi^\theta$ 
(smoothly depending on $\theta$), and define the solutions 
$(\tilde x^\theta, \tilde t^\theta, \tilde u^\theta, 
\tilde \alpha^\theta, 
\tilde \beta^\theta, \tilde p^\theta,\tilde q^\theta)(\Tilde X,
\Tilde Y)$
as in (\ref{TXY})--(\ref{TPQ}).   

\begin{figure}[htbp]
   \centering
\includegraphics[width=0.6\textwidth]{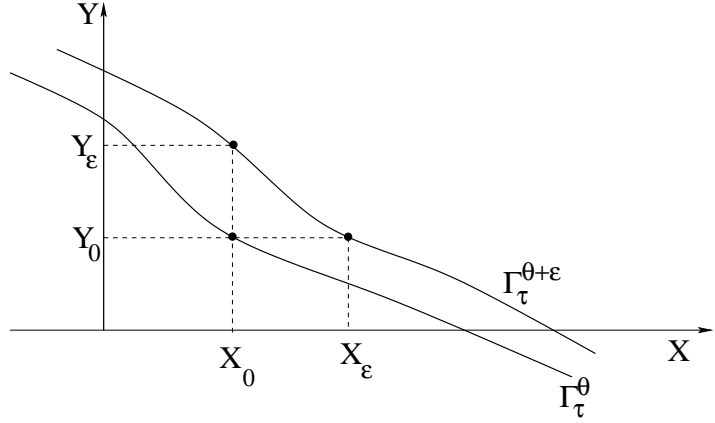}
   \caption{\small Given a representation 
   of the solutions $u^\theta$ in terms of the variables $X,Y$,
   the shifts $w,z$ are uniquely determined by (\ref{shifts}). 
   Here $\Gamma_\tau^\theta=\{(X,Y)\,;~~t^\theta(X,Y)=\tau\}$. 
    }
   \label{f:w28}
\end{figure}

 On the other hand, different relabelings of the $X,Y$ coordinates
determine different values for the integral in (\ref{ulg}).
Indeed, these correspond to different choices of the shifts $w,z$ 
in (\ref{fn1}).
To illustrate this point more clearly, fix a value  of $\theta$.
 Then, for $\ve>0$ small, 
the family of solutions
$u^{\theta+\ve}$ can be regarded as perturbations
of the solution $u^{\theta}$.    At a given point $(\tau,\bar x)$, the shifts 
$w(\tau,\bar x)$ and $z(\tau,\bar x)$ are uniquely determined as follows 
(Fig.~\ref{f:w28}).   Let $X_0, Y_0$ be the point in the $X$-$Y$ plane such that
$x^\theta(X_0, Y_0)=\bar x$, $t^\theta(X_0,Y_0) = \tau$.     For each $\ve>0$, 
define $X^\ve$ and $Y^\ve$ implicitly by setting
$$t^{\theta+\ve}(X_0, Y_\ve)~=~\tau\,,\qquad\qquad t^{\theta+\ve}(X_\ve, Y_0)~=~\tau
\,.$$
The shifts are then uniquely defined by setting
\bel{shifts}w(\tau,\bar x,)~=~\lim_{\ve\to 0} ~
{x^{\theta+\ve}(X_0, Y_\ve) - x^\theta(X_0, Y_0)\over\ve}
\,,\qquad
z(\tau,\bar x)~=~\lim_{\ve\to 0} ~{
x^{\theta+\ve}(X_\ve, Y_0) - x^\theta(X_0, Y_0)\over\ve}
\,.\eeq

The above considerations lead to
\v
{\bf Definition~3.} {\it The length $\|\gamma\|$ of the 
piecewise regular 
path $\gamma:\theta\mapsto (u_0^\theta, u_1^\theta)$ is defined as the infimum of the expressions in 
(\ref{ulg}), taken over all piecewise smooth relabelings 
of the $X$-$Y$ coordinates.}
\v
Based on the analysis in Section~3, we now
give an estimate on how the length of a regular path can grow in time.
\v
{\bf Theorem 5.} {\it Given any $K,T>0$, there exist
constants $\kappa_1,\ldots,\kappa_6$ in (\ref{ulg}) and $C_{K,T}>0$
such that the following holds.
Consider a path of solutions 
$\theta\mapsto (u^\theta, u_t^\theta)$  
of (\ref{2}), which is piecewise regular for $t\in [0,T]$ and 
where each $u^\theta$ has total energy $\leq K$.
Then its length satisfies the estimates
\bel{lengg}
\|\gamma^\tau\|~\leq~C_{K,T}\,
\|\gamma^0\|\qquad\qquad\hbox{for all}~~0\leq 
\tau\leq T\,.
\eeq}

{\bf Proof.} {\bf 1.}  To fix the ideas, let $u^\theta$ be structurally stable for every $\theta\in [0,1]\setminus\{\theta_1,\ldots,\theta_N\}$.

Fix $\ve>0$  and choose a 
relabeling of the variables $X,Y$ such that, at time $t=0$,
\bel{lee}
\int_0^1\left(\sum_{\ell=1}^6 \kappa_\ell \cdot
\int_{\Gamma_0^\theta}\Big\{ 
|J_\ell^\theta| \, \W^- dX~+ |H_\ell^\theta|\, \W^+\, dY\Big\}\right)
\, d\theta~\leq~\|\gamma^0\| + \ve\,.\eeq
Since the solution $u$ is smooth in the $X$-$Y$ 
variables and piecewise smooth in the $x$-$t$ variables,
the existence of the tangent vector is clear, for every 
$\theta\in [0,1]$ and $t\in [0,T]$.
We claim that, for every $\theta\notin \{\theta_1,\ldots,\theta_N\}$,
an estimate such as (\ref{es0}) holds. Namely
\bel{es00}\bega{l}\ds
 \Big\|(v^\theta(\tau), r^\theta(\tau), s^\theta(\tau))
 \Big\|_{(u^\theta(\tau),R^\theta(\tau),S^\theta(\tau))}\cr\cr
 \qquad \ds~\leq~
\exp\left\{C_0\tau+\int_0^\tau a^\theta(s)ds\right\}
\cdot\Big\|(v^\theta(0), r^\theta(0), s^\theta(0))
\Big\|_{(u^\theta(0),R^\theta(0),S^\theta(0))}\,.\enda\eeq
Here the constant $C_0$ and the integral of 
$a^\theta$ depend only on $T$
and on an upper bound on the total energy.

Integrating (\ref{es00}) over the interval $\theta\in[0,1]$,
one obtains an estimate of the form
$$
\|\gamma^\tau\|~\leq~C\,\bigl(\|\gamma^0\|+\ve\bigr)
\qquad\qquad\hbox{for all}~~0\leq 
\tau\leq T\,.
$$
This proves (\ref{lengg}), because  $\ve>0$ was arbitrary.
\v
{\bf 2.} It now remains to prove the estimate (\ref{es00}).
We observe that, if $u^\theta$ were smooth for all $(x,t)\in\R\times
[0,\tau]$, the result follows directly from (\ref{wn1}),
proved by the computations in Section~4.   
We need to show that the same conclusion can be reached
if $u^\theta$ is piecewise smooth, with structurally stable
singularities.

\begin{figure}[htbp]
   \centering
\includegraphics[width=0.4\textwidth]{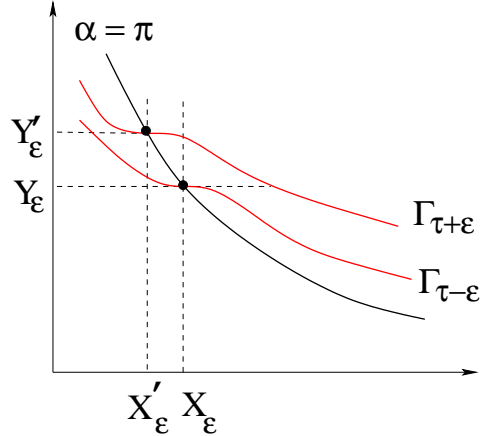}
   \caption{\small  Proving that the 
   rate of change in the length of a tangent vector is not affected 
   by the presence of a singularity.  
    }
   \label{f:w29}
\end{figure}

\begin{figure}[htbp]
   \centering
\includegraphics[width=0.6\textwidth]{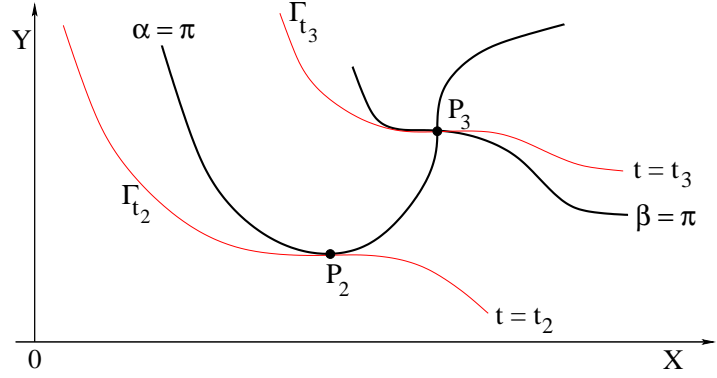}
   \caption{\small  Here $P_2$ is a singularity point of Type~2,  where $\alpha=\pi$ and $\alpha_X=0$,
   but $\alpha_{XX}\not= 0$ and $\beta\not= \pi$.
   At $P_3$ the solution has a singularity of Type~3, where
   $\alpha=\beta=\pi$, but $\alpha_X\not= 0$ and $\beta_Y\not= 0$.
    The weighted norm of the tangent vector is
    continuous at the times $t=t_2$ and $t=t_3$.}
   \label{f:w30}
\end{figure}

Fix a time $\tau$ and call   $\Gamma_\tau\doteq \{t^\theta(X,Y)=\tau\}$ 
the level set in the $X$-$Y$ plane.
Since the estimates of the previous section hold
in regions where $u^\theta$ is smooth, to obtain a bound on the 
weighted norm of the tangent vector it suffices to show that
the effect of isolated singularities is negligible.
To lighten the notation, in the following the superscript $^\theta$
will be omitted. 

With reference to Fig.~\ref{f:w29}, assume that the solution has a structurally 
stable singularity along a backward characteristic.    
We claim that this singularity does not affect the estimate
(\ref{wn1}).  In other words, the time derivative
$${d\over dt} \sum_{\ell=1}^6 \kappa_\ell \cdot\int_{\Gamma_t}
\Big\{|J_\ell| \, \W^-\,dX + |H_\ell|\, \W^+\, dY\Big\}$$
is not affected by the presence of the singularity.

For a given  time $\tau$, let 
$(X_\ve, Y_\ve)$ be the point where the curve $\Gamma_{\tau-\ve}=\{t(X,Y) = \tau-\ve\}$ 
intersects the singular curve  $\{\alpha(X,Y)=\pi\}$.
Similarly, let
$(X'_\ve, Y'_\ve)$ be the point where the curve 
$\Gamma_{\tau+\ve}=\{t(X,Y) = \tau+\ve\}$ 
intersects the singular curve  $\{\alpha(X,Y)=\pi\}$.

Define the curves
$$\left\{ \bega{rl}\sigma_\ve^+&\doteq~\Gamma_{\tau+\ve}\cap\{ X
\in [X_\ve', X_\ve]\}\,,\cr
\sigma_\ve^-&\doteq~\Gamma_{\tau-\ve}\cap\{ x\in [X_\ve', X_\ve]\}\,,\enda\right.
\qquad\qquad \left\{
\bega{rl}\eta_\ve^+&\doteq~\Gamma_{\tau+\ve}\cap\{ Y\in [Y_\ve, Y_\ve]\}\,,\cr
\eta_\ve^-&\doteq~\Gamma_{\tau-\ve}\cap\{ Y\in [Y_\ve, Y_\ve]\}\,.
\enda\right.$$
To prove our claim, it suffices  to show that
\bel{LX}
\lim_{\ve\to 0}~ {1\over\ve}\left( \int_{\sigma_\ve^+}- \int_{\sigma_\ve^-}\right)
\sum_{\ell=1}^6 |J_\ell| \W^-\, dX~=~0\,,\eeq
\bel{LY}
\lim_{\ve\to 0} ~{1\over\ve}\left( \int_{\eta_\ve^+}- \int_{\eta_\ve^-}\right)
\sum_{\ell=1}^6 |H_\ell| \W^+\, dY~=~0\,.\eeq

The first limit holds because the integrand is a continuous function of $X,Y$ and
$|X_\ve-X'_\ve|~=~\O(\ve)$.     
The second limit holds because the integrand is a continuous function of $X,Y$ and
$|Y_\ve'-Y_\ve|~=~\O(\ve)$. 
The basic estimate (\ref{wn1}) thus remains valid
also in the presence of singular curves  where $\alpha=\pi$
or $\beta=\pi$.
\v
Finally, we analyze what happens in the presence of singular points  
of Type 2, where $\alpha=\pi$ and $\alpha_x=0$, and of Type 3, where
$\alpha=\beta=\pi$.  Since the solution $u^\theta$
is structurally stable, there can be at most finitely many 
such points, say 
$$Q_j~=~(X_j, Y_j),\qquad\qquad j=1,\ldots,N\,.$$
To complete the proof of our claim, it thus suffices to show that,
at each time $\tau_j= t(X_j, Y_j)$, the map
\bel{ple}
t~\mapsto~
\int_0^1\left(\sum_{\ell=1}^6 \kappa_\ell \cdot
\int_{\Gamma_t}\Big\{ 
|J_\ell^\theta| \, \W^- dX~+ |H_\ell^\theta|\, \W^+\, dY\Big\}\right)
\eeq
is continuous at $t=\tau_j$.   But this is clear, because the path $\Gamma_t$ depends
continuously on $t$ and the integrands $J_\ell,H_\ell$ are uniformly bounded. Moreover, they are
continuous everywhere with a possible exception of
the finitely many singular points $Q_j$.
\endproof

\v
\v
\section{Construction of the geodesic distance}
\setcounter{equation}{0}

A key result proved in \cite{BChen} shows that
every path of solutions to (\ref{2}) can be approximated by a  path
which remains regular for $t\in [0,T]$.  More precisely, 
an application
of Thom's transversality theorem yields
\v
{\bf Theorem 6.} 
{\it Let the wave speed $c(u)$ satisfy the assumptions {\bf (A)}.
Let $(u^\theta,\alpha^\theta, \beta^\theta, p^\theta, q^\theta, 
x^\theta, t^\theta)(X,Y)$ be a path of $\C^\infty$ solutions
to the semilinear system (\ref{u})--(\ref{t}), depending smoothly on 
$\theta\in [0,1]$.
Then, for any  $T,\ve>0$ and any integer $k\ge 1$, 
there exists a perturbed  path 
of solutions $(\tilde u^\theta,\tilde \alpha^\theta, \tilde \beta^\theta, \tilde p^\theta, \tilde q^\theta, \tilde x^\theta, 
\tilde t^\theta)(X,Y)$ such that
\bel{approx}\Big\|(u^\theta-\tilde u^\theta,~\alpha^\theta-\tilde \alpha^\theta, ~\beta^\theta-\tilde \beta^\theta,~p^\theta- \tilde p^\theta, ~q^\theta-\tilde q^\theta, ~
x^\theta-\tilde x^\theta, ~t^\theta-\tilde t^\theta)\Big\|_{\C^k(\Omega)}
~<~\ve.
\eeq
Here $\Omega \subset\R^2$ is a domain containing the set
$$\Big\{(X,Y)\,;~~t^\theta(X,Y)\in [0,T]~~\hbox{or}~~
\tilde t^\theta(X,Y)\in [0,T]\,,
\qquad\hbox{for some}~ \theta\in [0,1]\Big\}.$$
Moreover, all except finitely many solutions $(\tilde u^\theta,\tilde \alpha^\theta, \tilde \beta^\theta, \tilde p^\theta, \tilde q^\theta, \tilde x^\theta, 
\tilde t^\theta)$ have structurally stable singularities 
inside $\Omega$.
}

In other words, by slightly perturbing the initial data 
$(u_0^\theta, u^\theta_1)$, $\theta\in [0,1]$, we can construct a 
one-parameter family of conservative solutions $u^\theta=u^\theta(t,x)$
which have structurally stable singularities, for all but finitely many
values of $\theta$. This implies 
that for all $t\in [0,T]$
the length of the path 
$\theta\mapsto u^\theta(t,\cdot)$
is well defined by the formula
\bel{length}
\|\gamma^t\|~\doteq~\int_0^1 \left\| {d\over d\theta} u^\theta(t)\right\|_{u^\theta(t)}
\, d\theta\,.\eeq 
Here $\|\cdot\|_u$ is a weighted norm defined as in 
(\ref{fn1})--(\ref{main}), or equivalently at (\ref{N2}).

\v

A  geodesic distance
$d^*$
on the space $H^1(\R)\times \L^2(\R)$ will be constructed in two steps.
\begi
\item[(i)] As proved in \cite{BChen}, there is an open dense set of initial data
\bel{Ddef}
\D~\subset~ \Big(\C^3(\R)\cap H^1(\R)\Big) \times\Big(\C^2(\R)\cap\L^2(\R)\Big),\eeq
such that, if $(u_0,u_1)\in\D$, then   the  solution of 
(\ref{0.1})-(\ref{0.2}) has 
structurally stable singularities.
On $\D^\infty\doteq \C^\infty_c\cap\D$ we  
construct a geodesic distance, defined as the infimum
among the weighted lengths of all piecewise regular paths 
connecting two given points.

\item[(ii)] By continuity, 
this distance can then  be extended  from  $\D^\infty$ to  
a larger space, defined as the completion of 
$\D^\infty$ w.r.t.~the distance $d^*$. In particular, this 
completion will contain the space $(H^1\cap W^{1,1})
\times (\L^2\cap\L^1)$.
\endi

More in detail, assume $(u_0, u_1)$, 
$(\tilde u_0, \tilde u_1)\in\D^\infty$.
Their total energies will be denoted by 
$$\E(u_0,u_1)~\doteq~\int \bigl[u_1^2+ c^2(u_0) u_{0,x}^2\bigr]\, dx\,,
\qquad\qquad \E(\tilde u_0,\tilde u_1)~\doteq~\int \bigl[\tilde u_1^2+ c^2(\tilde u_0) \tilde u_{0,x}^2\bigr]\, dx\,,$$
respectively.
Fix any constant $K>0$ and consider the subset of all 
data with energy $\leq K$, namely
\bel{XE}
X_{K}~\doteq~\Big\{(u_0,u_1)\in H^1(\R)\times \L^2(\R)\,;
\quad \E(u_0,u_1)\leq K\Big\}\,.\eeq
Notice that $X_K$ is positively invariant  for the flow 
generated by the wave equation.
\v
{\bf Definition 4.} {\it 
On $\D^\infty\cap X_K$ we  define the {\bf geodesic distance} 
$d^*\bigl((u_0, u_1), ~(\tilde u_0, \tilde u_1)\bigr)$  
as the infimum among
all weighted lengths 
of piecewise regular paths, which connect 
$(u_0, u_1)$
with  
$(\tilde u_0, \tilde u_1)$, always remaining inside $X_K$. 
Namely,
\bel{d*def}\bega{l}\ds
d^*\bigl((u_0, u_1), ~(\tilde u_0, \tilde u_1)\bigr)~\doteq~
\inf~\Big\{ \|\gamma\|\,;~~\gamma ~\hbox{is a piecewise
regular path}\,,\cr\cr
\qquad\qquad \gamma(0) = (u_0, u_1)\,,\quad
\gamma(1) = (\tilde u_0, \tilde u_1)\qquad 
\E(u_0^\theta, u_1^\theta)\leq K\quad\hbox{for all}~
\theta\in[0,1]\Big\}.\enda
\eeq
}

Since the concatenation of two piecewise regular paths is 
still a piecewise regular path (after a suitable re-parameterization),
it is clear that $d^*(\cdot,\cdot)$ is indeed a distance.
As a consequence of Theorem 5, we have
\v
{\bf Theorem 7.}   {\it Let the wave speed $c(\cdot)$ be smooth and satisfy (\ref{c0}).   Then the geodesic distance
$d^*$ renders Lipschitz continuous the flow generated by the
wave equation (\ref{0.1}).   In particular, let 
$(u_0,u_1)$ and $(\tilde u_0, \tilde u_1)$ be two initial data
in (\ref{0.2}). Then for all $t\in [0,T]$ 
the corresponding solutions satisfy
\bel{LP7}
d^*\Big( \bigl(u(t,\cdot), \, u_t(t,\cdot)\bigr),~\bigl(
\tilde u(t,\cdot), \, \tilde u_t(t,\cdot)\bigr)\Big)
~\leq~C_{K,T}\cdot 
d^*\Big((u_0,u_1),\,(\tilde u_0, \tilde u_1)\Big)\,.\eeq
Here $C_{K,T}$ is a constant depending only on $T$ and 
on an upper bound $K$ on the total energy.}
\v
{\bf Proof.} If the wave speed $c(\cdot)$ satisfies the generic
assumption ${\bf (A)}$ at (\ref{morse}), then the result is a direct consequence of 
Theorem~5.  To cover the general case,  it suffices to approximate
$c(\cdot)$ with 
a sequence of functions $c_n(\cdot)$
that satisfy the assumption {\bf (A)}. 
If $\|c_n-c\|_{C^3(\Omega)}~\to ~0$ as $n\to\infty$ for every bounded interval
$\Omega\subset\R$, then 
the flow generated by the velocities $c_n(\cdot)$
and the corresponding geodesic distances converge to the 
ones for $c(\cdot)$.
\endproof
\v
In the remainder of this section we compare the distance $d^*$ 
with more familiar distances in Sobolev spaces, and with a Wasserstein distance between energy measures.
\v
{\bf Proposition 2.}   {\it  There exists a constant 
$C'_K$ such that, for any $(u_0, u_1), \,(\tilde u_0, \tilde u_1)\in
\D^\infty\cap X_K$, 
\bel{d*w}
d^*\bigl((u_0, u_1), ~(\tilde u_0, \tilde u_1)\bigr)~
\leq~C'_K\cdot \Big( \|u_0-\tilde u_0\|_{H^1} + 
\|u_0-\tilde u_0\|_{W^{1.1}}+
\|u_1-\tilde u_1\|_{\L^2}+ \|u_1- \tilde u_1\|_{\L^1} \Big).\eeq
}

{\bf Proof.}
{\bf 1.} Define the function
\bel{Psdef}
\Psi(u)~\doteq~\int_0^u{c(s)}\, ds\,.\eeq
Observe that $\Psi:\R\mapsto \R$ is a smooth strictly increasing function, with smooth inverse $\Psi^{-1}$.
The total energy  can then be expressed as
$$\E(u_0,u_1)~\doteq ~\int \bigl[u_1^2 + c^2(u_0) u_{0,x}^2\bigr]\, dx
~=~\int \bigl[u_1^2 + \bigl(\Psi(u_0)_x\bigr)^2\bigr]\, dx\,.$$
Let $(\tilde u_0, \tilde u_1)$ be another initial data, with total energy $\Tilde\E$.
For $\theta\in [0,1]$, consider the interpolated data
$(u_0^\theta, u_1^\theta)$ where
\bel{interp}\left\{
\bega{rl}u_0^\theta&=~\Psi^{-1}\Big(\theta \Psi(\tilde u_0) + (1-\theta)\Psi(u_0)\Big),
\\[4mm]
u_1^\theta&=~\theta \tilde u_1 + (1-\theta) u_1\,.\enda\right.\eeq
When $\theta= 0,1$, it is clear that $(u_0^\theta, u_1^\theta)$
coincides with $(u_0, u_1)$ and $(\tilde u_0, \tilde u_1)$, respectively.
We check that the energy remains $\leq M$. 
Indeed,
\bel{EM}\bega{l}\ds
\int \bigl[(u_1^\theta)^2 + c^2(u_0^\theta) (u^\theta_{0,x})^2\bigr]\, dx
~=~\int \bigl[(u^\theta_1)^2 + \bigl(\Psi(u_0^\theta)_x\bigr)^2\bigr]\, dx\\[4mm]
\ds
\qquad =~\int \bigl[(\theta \tilde u_1 + (1-\theta) u_1)^2\, dx +\int
\bigl[ \theta \Psi(\tilde u_0)_x+ (1-\theta)  \Psi( u_0)_x\bigr]^2
\, dx\\[4mm]
\qquad \leq~\max\,
\bigl\{\E(u_0, u_1), \Tilde \E(\tilde u_0, \tilde u_1)\bigl\}~\leq~M.
\enda\eeq
\v
{\bf 2.} Next, we  estimate the weighted length of the path
$\gamma: \theta\mapsto (u_0^\theta, u_1^\theta)$ in (\ref{interp}), showing that
\bel{gle}
\|\gamma\|~\leq~ C\cdot  \Big( \|u_0-\tilde u_0\|_{H^1} + 
\|u_0-\tilde u_0\|_{W^{1.1}}+
\|u_1-\tilde u_1\|_{\L^2}+ \|u_1- \tilde u_1\|_{\L^1} \Big)\,,\eeq
for some constant $C$ depending only on the total energy.
To establish an upper bound for the weighted 
length $\|\gamma\|$, in the definition (\ref{main})
we choose the shifts $w=z=0$.
   In this way,
the integrals $I_1, I_4$, and $I_5$ vanish.  

We first calculate 
$(v^\theta, r^\theta,s^\theta)=
\frac{d}{d\theta}(u^\theta, R^\theta,S^\theta)$.  
From (\ref{Psdef}) it follows
\bel{Psi_pr}
\Psi'(u)~=~c(u)\,,
\eeq
\beq\label{Psi_inv}
(\Psi^{-1}(a))'~=~\frac{1}{\Psi'\big(\Psi^{-1}(a)\big)}
~=~\frac{1}{c\big(\Psi^{-1}(a)\big)}\,.
\eeq
Using (\ref{interp}) and (\ref{Psi_pr})-(\ref{Psi_inv}) we find
$$
v^\theta~=~\frac{d}{d\theta}u^\theta~=~
\frac{\Psi(\tilde u_0) -\Psi(u_0)}{c\big(\theta \Psi(\tilde u_0) + (1-\theta)\Psi(u_0)\big)}\,.$$
Since the wave speed $c(\cdot)$ is uniformly positive, the above implies
\bel{720}
{1\over K_1}\,|\tilde u_0 -u_0| ~\leq~|v^\theta|~
\leq~ K_1\,|\tilde u_0 -u_0|,
\eeq
for  a suitable constant $K_1$, depending on the function $c(\cdot)$
and on an upper bound for the energy.  
 
Next, we have
\bel{721}
R^\theta~=~u_1^\theta+\Psi(u_0^\theta)_x
~=~\theta \big({\tilde u}_1+ \Psi(\tilde u_0)_x\big) + 
(1-\theta) \big( u_1+ \Psi(u_0)_x\big)~=~\theta\tilde 
{R}+(1-\theta){R}.
\eeq
Hence
\bel{722}
r^\theta~=~\frac{d}{d\theta}R^\theta~=~
\big({\tilde u}_1+ \Psi(\tilde u_0)_x\big) 
-\big( u_1+ \Psi(u_0)_x\big)~=~\Tilde R-R\,.
\eeq
Similarly,
\bel{723}
s^\theta~=~\frac{d}{d\theta}S^\theta~=~
\big({\tilde u}_1- \Psi(\tilde u_0)_x\big) -\big( u_1- \Psi(u_0)_x\big)~=~\Tilde S-S.
\eeq
For later use, we observe that
\bel{724}
\int_0^1\Big(\int    |2R^\theta r^\theta|\,dx\, \Big)d\theta~=~
\int    |R-\tilde{R}|
\cdot\Big(\int_0^1 2|\tilde{R}^\theta|\,d\theta\Big)\,dx
\leq \int   |R-\tilde{R}|\cdot \bigl( |R|+|\tilde{R}|
\bigr)\,dx\,.
\eeq
Observing that the weights $\W^{\pm}$ satisfy a uniform bound
depending only on the total energy, and using (\ref{720})-(\ref{724}),
 we finally obtain
\bel{normes}\bega{rl}
\|\gamma\|
&=
~\ds
\int_0^1
\Big\|(v^\theta, r^\theta, s^\theta) \Big\|_{(u^\theta,R^\theta,S^\theta)}\, d\theta\cr\cr
&=\ds~
\int_0^1 
\Bigg\{\kappa_2\int   \Big\{   |r^\theta|(\W^-)^\theta  + 
 |s^\theta|(\W^+)^\theta    \Big\}dx\cr\cr
&\qquad \ds
\qquad+\kappa_3\int  |v^\theta|\,\Big\{\bigl(1+(R)^\theta)^2
\bigr)\,(\W^-)^\theta
+
\bigl(1+(S^\theta)^2\bigr)\,(\W^+)^\theta \Big\}\,dx 
\cr\cr
&\qquad\ds\qquad+\kappa_6\int   \Big\{ |2R^\theta r^\theta
|(\W^-)^\theta
+   |2S^\theta s^\theta|(\W^+)^\theta 
\Big\}\,dx\Bigg\}\, d\theta\cr\cr
&\leq
\ds~
K_2\cdot 
\bigg\{\int   \Big\{   |\tilde R-R| + 
 |\tilde S-S|    \Big\}dx +\|u_0-\tilde u_0\|_{\L^1}
 \cr\cr
&\quad \ds\qquad
 +\|u_0-\tilde u_0\|_{\L^\infty}\cdot\int_0^1\Big(\int  
 \Big\{(R^\theta)^2
+
(S^\theta)^2 \Big\}\,dx\Big)d\theta
\cr\cr
&\quad \ds\qquad
+\int   \Big\{ |R-\tilde{R}|\cdot
( |R|+|\tilde{R}|)
+   |S-\tilde{S}|\cdot( |S|+|\tilde{S}|)
\Big\}dx\bigg\}
\cr\cr
&\leq
\ds~
K_3\cdot
( \|u_0-\tilde u_0\|_{H^1} + 
\|u_0-\tilde u_0\|_{W^{1.1}}+
\|u_1-\tilde u_1\|_{\L^2}+ \|u_1- \tilde u_1\|_{\L^1} ),
\enda\eeq
where $K_2$ and $K_3$ are positive constants, depending on the 
upper bound for the energy.
In the last step, we use similar estimates as  in \eqref{EM}.
This completes the proof.
\endproof

\v
We conclude the paper by showing that the geodesic distance
$d^*$ in (\ref{d1}) controls both the $\L^1$ distance $\|u_0-\tilde u_0\|_{\L^1}$ and the Wasserstein distance between the corresponding energy measures $\mu,\tilde \mu$. 
\v
{\bf  Proposition 3.}
{\it There exists a constant $\delta_0$, depending only on an upper bound on the energy, such that for any $u_0, \tilde u_0\in H^1\cap \L^1$
and any $u_1,\tilde u_1\in \L^2$, one has
\bel{L1b}
\|u_0-\tilde u_0\|_{\L^1}~\leq ~\delta_0\cdot d^*\bigl( (u_0,u_1),\, 
(\tilde u_0, \tilde u_1)\bigr),\eeq
\bel{WD}
\sup_{\|f\|_{\C^1}\leq 1}\,\bigg|\int f\,d\mu - \int f d\tilde\mu\bigg|
~\leq~
\delta_0\cdot
d^*\bigl( (u_0,u_1),\, (\tilde u_0, \tilde u_1)\bigr)\,.\eeq
Here $\mu,\tilde \mu$  are the 
measures with densities
$u_1^2 + c^2(u_0) u_{0,x}^2$ and 
$\tilde u_1^2 + c^2(\tilde u_0) \tilde u_{0,x}^2$
w.r.t.~Lebesgue measure.}
\v
{\bf Proof.} {\bf 1.}
To prove  (\ref{L1b}) we first observe that
\bel{70}
|v|~\leq~ \left|v +\frac{Rw}{2c} - \frac{Sz}{2c}\right|+\left|\frac{Rw}{2c}\right| +\left|\frac{Sz}{2c}\right|
~\leq~ \left|v +\frac{Rw}{2c} - \frac{Sz}{2c}\right|+\frac{1}{4c}
\bigl|w(1+R^2)\bigr|+\frac{1}{4c}\bigl|z(1+S^2)\bigr|.
\eeq
The right hand side of (\ref{70})
is bounded by the integrands in $I_1$ and $I_3$ in \eqref{main}.
Recalling the definition (\ref{d*def}),
by (\ref{720}) for some constant $c_4>0$ we thus have
\bel{728}
\bega{l}
d^*\bigl((u_0, u_1), ~(\tilde u_0, \tilde u_1)\bigr)
~\geq ~\ds
c_4\cdot \inf_{\gamma}\bigg\{\int_0^1
\int  |v^\theta|\, dx\,d\theta\bigg\}
\cr\cr\ds
\qquad \leq~
c_4\cdot \inf_{\gamma}
 \int_0^1\left\|\frac{du^\theta}{d\theta}\right\|_{\L^1}
\,d\theta
~=~c_4\,\|u_0 - \tilde u_0\|_{\L^1}\,.
\enda\eeq
\v
{\bf 2.} Next, consider any regular path 
$\gamma:\theta\mapsto (u_0^\theta, u_1^\theta)$ joining $(u_0, u_1)$ with $(\tilde u_0, \tilde u_1)$.
Call $\mu^\theta$ the measure having  density
$(u_1^\theta)^2+ c^2(u_0^\theta) (u_0^\theta)^2=(R^\theta)^2+(S^\theta)^2$ w.r.t.~Lebesgue measure.

Then, for
any function $f$ such that $\|f\|_{\C^1}\leq 1$, one has
\bel{729}\bega{l}\ds\left|{d\over d\theta}\int f\,d\mu^\theta\right|\cr\cr
\leq
\ds~ K_5\cdot \int |f'|\cdot\Big\{   |w|\bigl(1+R^2\bigr)  +  
|z|\bigl(1+S^2\bigr)    \Big\}dx\cr\cr
\ds\quad+K_5\cdot \int  |f|\cdot
 \bigg\{ \Big|2R(r +w R_x)+R^2w_x
+2S(s+z S_x)+S^2z_x)\Big|    
\bigg\}dx
\cr\cr
\leq
\ds~ K_5\cdot \int \Big\{   |w|\bigl(1+R^2\bigr)  +  
|z|\bigl(1+S^2\bigr)    \Big\}dx\cr\cr
\ds\quad+K_5\cdot\int  
 \bigg\{ \Big|2R(r +w R_x)+R^2w_x+\frac{c'}{4c^2}(R^2S-S^2R)(w-z)
\Big|
\cr\cr
\ds\quad\qquad\qquad
+   \Big|2S(s+zS_x)+S^2z_x+\frac{c'}{4c^2}(S^2R-R^2S)(w-z)\Big|
\bigg\}dx\,.
\enda
\eeq
Using \eqref{rstt}, we see that the two integrals on the right hand side of (\ref{729}) 
are exactly $I_1$ and $I_6$ without potential terms $\W^-$ and $\W^+$, hence are dominated by the integrals in (\ref{main}). 
Integrating w.r.t.~$\theta\in [0,1]$, one obtains
(\ref{WD}).
\endproof
\v
{\bf Acknowledgment.} This research  was partially supported
by NSF, with grant  DMS-1411786: ``Hyperbolic Conservation Laws and Applications".

\end{document}